\RequirePackage{etoolbox}
\csdef{input@path}{{style/}{graphics/}}
\documentclass[aos,MSNbibl,nameyear,seceqn,dvips]{arximspdf}
\usepackage{graphicx}
\doi{10.1214/14-AOS1268}
\volume{43}
\issue{2}
\pubyear{2015}
\firstpage{785}
\lastpage{818}
\docsubty{FLA}

\makeatletter
\newcommand{\lleft}{\left}
\newcommand{\rright}{\right}
\newcommand{\rrVert}{\Vert}
\newcommand{\llVert}{\Vert}
\newtheorem{thmm}{Theorem}[section]
\newtheorem{lemma}{Lemma}[section]
\newtheorem{corollary}{Corollary}[section]
\newtheorem{proposition}{Proposition}[section]
\newproclaim{remark}{Remark}[section]
\makeatother

\begin{document}
\begin{frontmatter}

\title{Rate-optimal posterior contraction for sparse PCA}
\runtitle{Bayes sparse PCA}

\begin{aug}
\author[A]{\fnms{Chao} \snm{Gao}\ead[label=e1]{chao.gao@yale.edu}}
\and
\author[A]{\fnms{Harrison H.} \snm{Zhou}\corref{}\ead[label=e2]{huibin.zhou@yale.edu}}
\runauthor{C. Gao and H.~H. Zhou}
\affiliation{Yale University}
\address[A]{Department of Statistics\\
Yale University\\
New Haven, Connecticut 06511\\
USA\\
\printead{e1}\\
\phantom{E-mail:\ }\printead*{e2}}
\end{aug}

%
\received{\smonth{11} \syear{2013}}
%
\revised{\smonth{8} \syear{2014}}

%
\begin{abstract}
Principal component analysis (PCA) is possibly one of the most widely
used statistical tools to recover a low-rank structure of the data. In
the high-dimensional settings, the leading eigenvector of the sample
covariance can be nearly orthogonal to the true eigenvector. A sparse
structure is then commonly assumed along with a low rank structure.
Recently, minimax estimation rates of sparse PCA were established under
various interesting settings. On the other side, Bayesian methods are
becoming more and more popular in high-dimensional estimation, but
there is little work to connect frequentist properties and Bayesian
methodologies for high-dimensional data analysis. In this paper, we
propose a prior for the sparse PCA problem and analyze its theoretical
properties. The prior adapts to both sparsity and rank. The posterior
distribution is shown to contract to the truth at optimal minimax
rates. In addition, a computationally efficient strategy for the
rank-one case is discussed.
\end{abstract}
%
%
\end{frontmatter}

\section{Introduction}

Principal component analysis is a classical statistical tool used to
project data
into a lower dimensional space while maximizing the variance
[\citeauthor{jolliffe86} (\citeyear{jolliffe86})]. When the sample size
$n$ is
small compared to the number of variables $p$, \citeauthor{johnstone09}
(%
\citeyear{johnstone09}) show that the standard PCA may fail in the sense
that the leading eigenvector of the sample covariance can be nearly
orthogonal to the true eigenvector. Therefore, the recovery of principal
components in the high-dimensional setting requires extra structural
assumptions. The sparse PCA, assuming that the leading eigenvectors or
eigen-subspace only depend on a relatively small number of variables, is
applied in a wide range of applications. Estimation methods for sparse PCA
problems are proposed in \citeauthor{zou06} (\citeyear{zou06}) and %
\citeauthor{d07} (\citeyear{d07}). \citeauthor{amini09} (\citeyear{amini09})
and \citeauthor{ma13} (\citeyear{ma13}) obtain rates of convergence of
sparse PCA methods under the spiked covariance model proposed in %
\citeauthor{johnstone09} (\citeyear{johnstone09}). Minimax rates of sparse
PCA problems are established by \citeauthor{birnbaum13} (%
\citeyear{birnbaum13}), \citeauthor{cai12} (\citeyear{cai13,cai12}) and \citeauthor{vu12}
(\citeyear
{vu12}%
) under various interesting settings.

Bayesian methods have been very popular in high-dimensional estimation, but
there is little work to connect frequentist properties and Bayesian
methodologies for high-dimensional models. This paper serves as a bridge
between the frequentist and Bayesian worlds by addressing the following
question for high-dimensional PCA: Is it possible for a Bayes procedure to
optimally recover the leading principal components in the sense that the
posterior distribution contracts to the truth with a minimax rate? The
optimal posterior contraction rate immediately implies that the posterior
mean attains the optimal convergence rate as a point estimator.

In this paper we consider a spiked covariance model with an unknown growing
rank. We propose a sparse prior on
the covariance matrix with a spiked structure and show that the induced
posterior distribution contracts to the truth with an optimal minimax rate.
The assumptions are nearly identical to those in \citeauthor{vu12} (%
\citeyear{vu12}), where the rank of the principal space $r=O(\log p)$ and
the number of nonzero entries of each spike $s$ is allowed to be at the
order of $p^{1-c}$ for any $c\in(0,1)$, as long as the minimax rate
$\frac{%
rs\log p}{n}\rightarrow0$. In addition, we prove that the posterior
distribution consistently estimates the rank. To the best of our knowledge,
this is the first work where a Bayes procedure is able to adapt to both the
sparsity and the rank.

There are two key ingredients in our approach. The first ingredient is in
the design of the prior. We propose a prior that imposes a spiked structure
on a random covariance matrix, under which each spike is sparse and
orthogonal to each other. This leads to sufficient prior concentration
together with the sparse property. In addition, each spike has a
bounded $%
l^{2}$ norm under the prior distribution such that there is a fixed
eigen-gap between the spikes and the noise, which eventually leads to
consistent rank estimation. The second ingredient is in constructing
appropriate tests in the proof of posterior contraction under spectral and
Frobenius norms. We first construct a test with the alternative hypothesis
outside of the neighborhood of the true covariance under the spectral norm.
For the covariance matrices inside the neighborhood of the truth under the
spectral norm, we propose a delicate way to divide the region into many
small pieces, where the likelihood ratio test is applicable in each small
region. A final test is then constructed by combining these small
tests. The
errors are controlled by correctly calculating the covering number
under the
metric for measuring the distance of subspaces.

The theoretical tools we use for this problem follow the recent line of
developments in Bayesian nonparametrics pioneered by \citeauthor
{barron88} (%
\citeyear{barron88}) and \citeauthor{barron99} (\citeyear
{barron99}), which
generalize the testing theory of \citeauthor{lecam73} (\citeyear{lecam73})
and \citeauthor{schwartz65} (\citeyear{schwartz65}) to construct an
exponentially consistent test on the essential support of a prior to prove
posterior consistency. The idea was later extended by \citeauthor
{ghosal00} (%
\citeyear{ghosal00}) and \citeauthor{shen01} (\citeyear{shen01}) to prove
rates of convergence of posterior distribution. Compared to Bayesian
nonparametrics, little work has been done for Bayesian high-dimensional
estimation, especially in the sparse setting. \citeauthor{castillo12} (
\citeyear{castillo12}) is the first work in this area. They prove
rates of
convergence in sparse vector estimation for a large class of priors.

The works closely related to this paper are \citeauthor{banerjee13} (%
\citeyear{banerjee13}) and \citeauthor{pati12} (\citeyear{pati12}). %
\citeauthor{banerjee13} (\citeyear{banerjee13}) study rates of convergence
for Bayesian precision matrix estimation by considering a conjugate prior.
But as discussed in \citeauthor{birnbaum13} (\citeyear{birnbaum13}),
estimation of sparse or bandable covariance/precision matrix is different
from that of sparse principal subspace. The optimal rates of
convergence can
be different. \citeauthor{pati12} (\citeyear{pati12}) study Bayesian
covariance matrix estimation for a sparse factor model, which is similar
to the spiked covariance model in the PCA problem. Instead of
estimating the principal subspace as in the
PCA problem, they consider estimating the whole covariance matrix. The
posterior rate of convergence they obtain is not optimal, especially
when the rank $r$ is allowed to grow with the sample size $n$.

The paper is organized as follows. In Section~\ref{sec:setting}, we
introduce the sparse PCA problem and define the parameter space. In
Section~\ref{sec:main}, we propose a prior and state the main result of the
posterior convergence. Section~\ref{sec:disc} introduces an algorithm to
compute the posterior mean in the rank-one case along with other
discussions. All the proofs are presented in Section~\ref{sec:proof}, with
some technical results given in the supplementary material [\citeauthor
{gaosupp} (\citeyear{gaosupp})].

\section{The sparse PCA}
\label{sec:setting}

Let $X_{1},\ldots,X_{n}$ be i.i.d. observations from $P_{\Sigma
}=N(0,\Sigma)$%
, with $\Sigma$ being a $p\times p$ covariance matrix with a spiked
structure
\[
\Sigma=\sum_{l=1}^{r}\theta_{l}
\theta_{l}^{T}+I_{p\times p},
\]
where $\theta_{l}^{T}\theta_{k}=0$ for any $l\neq k$. It is easy to see
that $ ( \Vert \theta_{1}\Vert ^{-1}\theta_{1},\ldots,\Vert \theta
_{r}\Vert ^{-1}\theta
_{r} ) $ are the first $r$ eigenvectors of $\Sigma$, with the
corresponding eigenvalues $ ( \Vert \theta_{1}\Vert ^{2}+1,\ldots,\Vert \theta
_{r}\Vert ^{2}+1 ) $. The rest $p-r$ eigenvalues are all $1$. The spiked
covariance is proposed by \citeauthor{johnstone09} (\citeyear{johnstone09})
to model data with a sparse and low-rank structure. An equivalent
representation of the data is
%
\begin{equation}
X_{i}=V_{0}\Lambda_{0}^{1/2}W_{i}+Z_{i}\qquad
\mbox{for }i=1,2,\ldots,n, \label{eq:latent}
\end{equation}
where $W_{i}\sim N(0,I_{r\times r})$ and $Z_{i}\sim N(0,I_{p\times p})$ are
independent. The matrix $V_{0}$ is defined as $V_{0}= [ \Vert \theta
_{1}\Vert ^{-1}\theta_{1},\ldots,\Vert \theta_{r}\Vert ^{-1}\theta_{r} ] $ and
$%
\Lambda_{0}=\operatorname{diag} ( \Vert \theta
_{1}\Vert ^{2},\ldots,\Vert \theta
_{r}\Vert ^{2} ) $. In such latent variable representation,
$V_{0}\Lambda
_{0}^{1/2}W_{i}$ models the signal part, which lives in an $r$-dimensional
subspace, and $Z_{i}$ is the noise part, which has the same variance on
every direction. Since the $r$-dimensional subspace is determined by its
projection matrix $V_{0}V_{0}^{T}$, the goal here is to recover the
principal subspace by estimating its projection matrix in the Frobenius
loss,
\[
\bigl\Vert \hat{V}\hat{V}^{T}-V_{0}V_{0}^{T}
\bigr\Vert _{F}.
\]

In a high-dimensional setting, extra structural assumptions are needed for
consistent estimation. We assume that the first $r$ eigenvectors are sparse,
in the sense that each of them only depends on a few coordinates among the
total number $p$. Define $S_{0,l}=\operatorname{supp}(\theta_{l})$ for
$l=1,2,\ldots,r$%
, the support of the $l$th eigenvector. We assume $l^0$ sparsity on
each spike by
$\max_{1\leq l\leq r}|S_{0,l}|\leq s$. The parameter space for the
covariance matrix is
\[
\mathcal{G}(p,s,r)=\lleft\{ %
\begin{array} {c} \displaystyle\Sigma=\sum
_{l=1}^{r}\theta_{l}
\theta_{l}^{T}+I\dvtx \max_{1\leq l\leq
r}|S_{0l}|
\leq s,\theta_{l}\in\mathbb{R}^{p},
\\
\displaystyle\theta_{l}^{T}\theta_{k}=0\mbox{ for }k\neq l,
\Vert \theta_{l}\Vert ^{2}\in \bigl(K^{-1},K
\bigr)%
\end{array} %
 \rright\},
\]
where $K>0$ is a constant, which we treat as being known in this paper.
The sparsity we consider matches the column sparsity in \citeauthor
{vu12} (\citeyear{vu12}) in the $l^0$ case. We require both upper and
lower bounds for $\Vert \theta_l\Vert ^2$. The lower bound implies an eigengap,
which leads to rank adaptation and subspace estimation, while the upper
bound controls the spectral norm of $\Sigma$, which leads to estimation
of the whole covariance matrix.
\citeauthor{vu12} (\citeyear{vu12}) prove that under the
assumptions
\[
r\leq m\log p\quad\mbox{and}\quad s\leq p^{1-c}\qquad\mbox{for some constants }c\in
(0,1)%
\mbox{ and }m>0,
\]
the minimax rate\footnote{The minimax rate is obtained by combining
Theorem~3.5 and Corollary~3.2 in \citeauthor{vu12} (\citeyear{vu12}).
The upper bound is a special case of their Corollary~3.2 because our
parameter space is a subset of theirs. The lower bound holds by
observing that the least favorable class in the proof of their Theorem~3.5 is a subset of our parameter space.} of principal subspace
estimation is
\[
\inf_{\hat{V}}\sup_{\Sigma\in\mathcal{G}(p,s,r)}P_{\Sigma
}^{n}
\bigl\Vert \hat {V}%
\hat{V}^{T}-V_{0}V_{0}^{T}
\bigr\Vert _{F}^{2}\asymp\frac{rs\log p}{n}.
\]
The goal of this paper is to prove an alternative result, adaptive Bayesian
estimation, by designing an appropriate prior $\Pi$, such that
%
\begin{equation}\qquad
\sup_{\Sigma\in\mathcal{G}(p,s,r)}P_{\Sigma}^{n}\Pi \bigl( \bigl\Vert
VV^{T}-V_{0}V_{0}^{T}\bigr\Vert
_{F}^{2}>M\varepsilon^{2}|X^{n}
\bigr) \leq \delta\qquad \mbox{for some }M>0, \label{eq:postconverge}
\end{equation}
where $\varepsilon^{2}=\frac{rs\log p}{n}$ is the minimax rate and
$X^{n}\sim
P_{\Sigma}^{n}$. The number $\delta>0$ satisfies $\lim_{(n,s,p,r)%
\rightarrow\infty}\delta=0$. The posterior contraction (\ref
{eq:postconverge}) leads to a risk bound of a point estimator. Let
$\mathbb{E}_{\Pi}$ be the expectation
under the prior distribution $\Pi$. Consider the posterior mean of the
subspace projection matrix $\mathbb{E}_{\Pi} (
VV^{T}|X^{n} ) $.
Its risk upper bound is given in the following proposition. We prove the
proposition in the supplementary material [\citeauthor {gaosupp}
(\citeyear {gaosupp})].

\begin{proposition}
\label{prop:pointestimate} Equation (\ref{eq:postconverge}) implies
\[
\sup_{\Sigma\in\mathcal{G}(p,s,r)}P_{\Sigma}^{n}\bigl\llVert
\mathbb{E} 
_{\Pi} \bigl( VV^{T}|X^{n}
\bigr) -V_{0}V_{0}^{T}\bigr\rrVert
_{F}^{2}\leq M\varepsilon^{2}+2(p+r)\delta.
\]
\end{proposition}

\begin{remark}
In this paper, the number $\delta$ in (\ref{eq:postconverge}) is at
an order
of $\exp ( -C^{\prime}n\varepsilon^{2} )$ for some
$C^{\prime
}>0$%
. Thus the dominating term in $M\varepsilon^{2}+2(p+r)\delta$ is
$M\varepsilon
^{2}$. The posterior mean is a rate-optimal point estimator.
\end{remark}

\begin{remark}
The matrix $\mathbb{E}_{\Pi} ( VV^{T}|X^{n} ) $ may not be a
projection matrix. However, it is still a valid estimator of the true
projection matrix $V_{0}V_{0}^{T}$. A~projection matrix estimator can be
obtained by projecting the posterior mean $\mathbb{E}_{\Pi} (
VV^{T}|X^{n} ) $ to the space of projection matrices under the
Frobenius norm. Denote the projection by $\hat{V}\hat{V}^{T}$. It can be
shown that $\Vert \hat{V}\hat{V}^{T}-V_{0}V_{0}^{T}\Vert _{F}\leq2\llVert
\mathbb{E}_{\Pi} ( VV^{T}|X^{n} ) -V_{0}V_{0}^{T}\rrVert _{F}$.
\end{remark}

\subsection{Notation}

In this paper, we use $\Gamma$ to denote a $p\times p$ spiked covariance
matrix with structure $\Gamma=AA^{T}+I$, where $A=[\eta_{1},\eta
_{2},\ldots,\eta_{\xi}]$ is a $p\times\xi$ matrix with orthogonal columns.
We use $S_{l}$ to denote the support of $\eta_{l}$ for each
$l=1,2,\ldots,\xi$%
. Define
\begin{eqnarray*}
V&=& \bigl[ \Vert \eta_{1}\Vert ^{-1}\eta_{1},
\Vert \eta_{2}\Vert ^{-1}\Vert \eta _{2}\Vert ,
\ldots,\Vert \eta_{\xi}\Vert ^{-1}\eta_{\xi} \bigr] ,\\
\Lambda &=&\operatorname{%
diag} \bigl( \Vert \eta_{1}\Vert
^{2},\Vert \eta_{2}\Vert ^{2},\ldots,\Vert
\eta_{\xi
}\Vert ^{2} \bigr).
\end{eqnarray*}
Then $V$ is a $p\times\xi$ unitary matrix, and $\Gamma$ has an
alternative representation $\Gamma=V\Lambda V^{T}+I$. We use
$P_{\Gamma}$
to denote the probability or the expectation under the multivariate normal
distribution $N(0,\Gamma)$ and $P_{\Gamma}^{n}$ to denote the product
measure. The symbol $\mathbb{P}$ stands for a generic probability whose
distribution will be made clear through the context. Correspondingly,
we use
$(\Sigma,A_{0},r,\theta_{l},S_{0l},V_{0},\Lambda_{0})$ to denote the true
version of $(\Gamma,A,\xi,\eta_{l},S_{l},V,\Lambda)$.

For a matrix $A$, we use $\Vert A\Vert $ to denote its spectral norm and $\Vert A\Vert _{F}$
for the Frobenius norm. We define $\mathcal{U}(d,r)$ to be the space of
all $%
d\times r$ unitary matrices for $d\geq r$ such that for any $U\in
\mathcal{U}(d,r)$, $U^TU=I_{r\times r}$. For any $U,V\in\mathcal
{U}(d,r)$%
, define the distance $d_{\Lambda}(\cdot,\cdot)$ by $d_{\Lambda
}(\cdot
,\cdot)=\Vert U\Lambda U^{T}-V\Lambda V^{T}\Vert _{F}$ for some diagonal
matrix $\Lambda$. We omit the subscript $%
\Lambda$ and write $d(\cdot,\cdot)=d_{\Lambda}(\cdot,\cdot)$
whenever $%
\Lambda=I$. The number $\varepsilon^{2}$ stands for the minimax rate
$\frac{%
rs\log p}{n}$ throughout the paper.

\section{The prior and the main results}
\label{sec:main}

\begin{figure}

\includegraphics{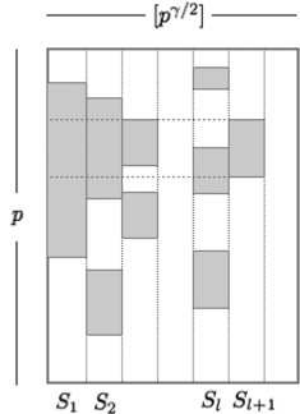}

\caption{An illustration of the prior. The shaded areas are $\{S_l\}
_{l=1}^{[p^{\gamma/2}]}$. The parts inside the
dashed lines correspond to $u_1,\ldots,u_l$ defined in (\protect\ref
{eq:dashed}).}
\label{fig:Prior}
\end{figure}

We propose a prior $\Pi$ from which we can sample a random covariance
matrix with structure $\Gamma=AA^{T}+I=\sum_{l=1}^{\xi}\eta_{l}\eta
_{l}^{T}+I$, where $A$ is a $p\times\xi$ matrix. The prior $\Pi$ is
described as follows:
\begin{longlist}[(1)]
\item[(1)] for each $l\in\{1,\ldots,[p^{\gamma/2}] \}$, we randomly
choose $%
S_{l}\subset\{1,\ldots,p\}$ by letting the indicator $\mathbb{I}\{i\in
S_{l}\}$
for each $i=1,\ldots,p$ follow a Bernoulli distribution with parameter $%
p^{-(1+\gamma)}$;
\item[(2)] given $(S_{1},\ldots,S_{[p^{\gamma/2}]})$, we sample a
$p\times
[p^{\gamma/2}]$ matrix $\bar{A}=[{\eta}_1,\ldots,{\eta}_{[p^{\gamma/2}]}]$
from $G_{(S_{1},\ldots,S_{[p^{\gamma/2}]})}$ to be specified below, and
then let
$\Gamma=\bar{A}\bar{A}^{T}+I$.
\end{longlist}
The $p\times[p^{\gamma/2}]$ matrix $\bar{A}$ (Figure \ref{fig:Prior}) may contain some zero
columns under the above sampling procedure. With slight abuse of
notation, we gather those nonzero columns to form the matrix $A=[\eta
_1,\ldots,\eta_{\xi}]$, with $S_l$ being the support of the column $\eta
_l$. Note that $\Gamma=\bar{A}\bar{A}^T+I=AA^T+I$, where $A$ is a
$p\times\xi$ matrix.
After specifying the distribution $G_{(S_{1},\ldots,S_{[p^{\gamma/2}]})}$,
the number of nonzero columns $\xi$ is also the rank of $A$.

\begin{remark}
The number $\gamma>0$ is a fixed constant in the prior. With $%
p^{-(1+\gamma)} $ as the mean for $\mathbb{I}\{i\in S_l\}$, the
cardinality $%
|S_l|$ is small with high probability under the prior distribution.
\end{remark}

\begin{remark}
The number $[p^{\gamma/2}]$ is an upper bound of the rank $\xi$. In
this paper, we assume that the true rank $r$ is at the order of $O(\log
p)$. Since $\log p\ll p^{[\gamma/2]}$, the range of $\xi$ covers the
range of $r$.
\end{remark}

We need to define a distribution $G_{d}^{\ast}$ on $\mathbb{R}^{d}$
to help
introduce $G_{(S_{1},\ldots,S_{[p^{\gamma/2}]})}$. Let
$Z=(Z_{1},\ldots,Z_{d})$ follow $N(0,I_{d\times d})$ and $U$ follow the
uniform distribution on the interval $%
[(2K)^{-1/2},(2K)^{1/2}]$. Then $G_{d}^{\ast}$ is defined to be the
distribution of
%
\begin{equation}
\biggl( \frac{UZ_{1}}{\Vert Z\Vert },\ldots,\frac{UZ_{d}}{\Vert Z\Vert } \biggr). \label{eq:priorcons}
\end{equation}
Now we are ready to specify the random matrix prior
$G_{(S_{1},\ldots,S_{[p^{\gamma/2}]})}$, which induces a distribution over
the matrix $\bar{A}=[\eta_{1},\eta
_{2},\ldots,\eta_{[p^{\gamma/2}]}]$. For any vector $v$ and any subset
$S$, we use the notation $v^T=(v_S^T,v_{S^c}^T)$. We describe the prior
through a sequential sampling
procedure. If $|S_1|=0$, we set $\eta_1=0$. Otherwise, we sample $\eta
_{1,S_{1}}\sim G_{|S_{1}|}^{\ast}$ and
let
\[
\eta_{1}=\pmatrix{\eta_{1,S_{1}}
\cr
0}.
\]
Suppose we have already obtained $(\eta_{1},\ldots,\eta_{l})$ and then
sample $%
\eta_{l+1}$, conditioning on $(\eta_{1},\ldots,\eta_{l})$. We set $\eta
_{l+1,S_{l+1}^{c}}=0$. The prior distribution of $\eta_{l+1,S_{l+1}}$
depends on $\eta_{i},1\leq i\leq l$, through values of $\eta_{i}$'s
on the
index set $S_{l+1}$. For simplicity, denote
%
\begin{equation}
(u_{1},\ldots,u_{l})=(\eta_{1,S_{l+1}},\ldots,
\eta_{l,S_{l+1}}). \label{eq:dashed}
\end{equation}
Define $l^{\ast}=\dim ( \operatorname{span} \{
u_{1},\ldots,u_{l} \}
 ) $. If $|S_{l+1}|-l^*\leq0$, we set $\eta_{l+1,S_{l+1}}=0$.
Otherwise, let $H_{l}$ be the projection matrix from $\mathbb{R}^{S_{l+1}}
$ to the subspace spanned by $ \{ u_{1},\ldots,u_{l} \} $. There
is a bijective
linear isometry $T_{l}$ induced by $H_{l}$ such that
\[
T_{l}\dvtx (I-H_{l})\mathbb{R}^{S_{l+1}}\rightarrow
\mathbb {R}^{|S_{l+1}|-l^{\ast
}},\qquad  T_{l}^{-1}\dvtx
\mathbb{R}^{|S_{l+1}|-l^{\ast}}\rightarrow (I-H_{l})%
\mathbb{R}^{S_{l+1}}.
\]
Remember that a linear isometry preserves the norms in the sense that
$\Vert T_lv\Vert =\Vert v\Vert $.
We sample $\bar{u}_{l+1}$ from $G_{|S_{l+1}|-l^{\ast}}^{\ast}$ and
let $%
u_{l+1}=T_{l}^{-1}\bar{u}_{l+1}$. Set $\eta_{l+1,S_{l+1}}=u_{l+1}$. Then
we have specified $\eta_{l+1}^{T}$, which is $(\eta
_{l+1,S_{l+1}}^{T},0^{T})$. Repeating this step, we obtain $A=[\eta
_{1},\ldots,\eta_{[p^{\gamma/2}]}]$. The prior $\Pi$ on the random
covariance matrix $%
\Gamma$ is now fully specified.

After collecting the nonzero $\eta_l$'s, we observe that the prior
$\Pi
$ explicitly samples a spiked covariance matrix $%
\Gamma=\bar{A}\bar{A}^T+I=AA^T+I =\sum_{l=1}^{\xi}\eta_{l}\eta
_{l}^{T}+I$ with the number of spikes
being $\xi$. The prior $\Pi$ imposes orthogonality on the spikes,
since $%
\eta_{l+1}$ is sampled on the orthogonal complement of the space
$\operatorname
{span%
} \{ \eta_{1},\eta_{2},\ldots,\eta_{l} \} $. Therefore,
$\eta
_{k}^{T}\eta_{l}=0$ for each $k\neq l$, and $\{\Vert \eta_{l}\Vert ^{-1}\eta
_{l}\}_{l=1}^{\xi}$ are the eigenvectors. For each eigenvector $\Vert \eta
_{l}\Vert ^{-1}\eta_{l}$, its support is in $S_{l}$, whose cardinality is small
under the prior distribution. Moreover, the first $\xi$ eigenvalues
are all
bounded from $1$ and $\infty$ because $\Vert \eta_{l}\Vert ^{2}\in{}[
(2K)^{-1},(2K)]$.

Given the data $X^{n}=(X_{1},\ldots,X_{n})\sim P_{\Sigma}^{n}$, the posterior
distribution is defined as
%
\begin{equation}
\Pi \bigl( B|X^{n} \bigr) ={\int_{B}\frac{dP_{\Gamma
}^{n}}{dP_{\Sigma
}^{n}}\bigl(X^{n}\bigr)\,d\Pi(\Gamma)}\Big/\biggl(\int\frac{dP_{\Gamma}^{n}}{dP_{\Sigma
}^{n}}%
\bigl(X^{n}\bigr)\,d\Pi(\Gamma)\biggr), \label{eq:postdef}
\end{equation}
for any measurable set $B$. The following theorem is the main result of this
paper. The posterior distribution contracts to the truth with an optimal
minimax rate.

\begin{thmm}
\label{thmm:postsub} Assume $\varepsilon\rightarrow0$, $r\leq m (
s\wedge
\log p ) $ and $n\leq p^{m}$ for some constant $m>0$. Then there exists
$M_{\gamma,K,m}^{\prime}>0$, such that for any $M^{\prime}>M_{\gamma
,K,m}^{\prime}$, we have
\[
\sup_{\Sigma\in\mathcal{G}(p,s,r)}P_{\Sigma}^{n}\Pi \bigl( \bigl\Vert
VV^{T}-V_{0}V_{0}^{T}\bigr\Vert
_{F}>M^{\prime}\varepsilon|X^{n} \bigr) \leq \exp
\bigl( -C_{(\gamma,K,m,M)}n\varepsilon^{2} \bigr) ,
\]
for some constant $C_{(\gamma,K,m,M^{\prime})}>0$ only depending on $%
(\gamma,K,m,M^{\prime})$.
\end{thmm}

Note that we have obtained the optimal posterior contraction rate under a
``mildly growing rank'' regime $r\leq m\log p$, which is also
assumed in \citeauthor{vu12} (\citeyear{vu12}), for them to match the upper
and lower bounds for minimax estimation. The assumption $n\leq p^{m}$ is
a convenient but mild condition in high-dimensional statistics to prove
rates of convergence in
expectation rather than with high probability; see \citeauthor {cai11}
(\citeyear
{cai11}), \citeauthor {paul12} (\citeyear {paul12}), etc. The
posterior contraction
result implies the same rate of convergence in expectation of a point
estimator (Corollary~\ref{cor:pointconverge}), and thus we need such an
assumption to hold. Additionally, we
assume $r\leq ms$, which means that the level of the rank is not above the
level of sparsity. This assumption is due to the fact that $V_{0}$ can be
only identified up to a unitary transformation, that is, $%
V_{0}V_{0}^{T}=(V_{0}Q)(V_{0}Q)^{T}$ for any $Q\in\mathcal{U}(r,r)$, and
for some $Q$ such that each row of $V_{0}Q$ may have at least $r$
nonzero entries.

As shown in Proposition~\ref{prop:pointestimate}, we can use the posterior
mean as a point estimator to achieve the minimax optimal rate of convergence.

\begin{corollary} \label{cor:pointconverge}
Under the setting of Theorem~\ref{thmm:postsub}, we have
\[
\sup_{\Sigma\in\mathcal{G}(p,s,r)}P_{\Sigma}^{n}\bigl\llVert
\mathbb{E} 
_{\Pi} \bigl( VV^{T}|X^{n}
\bigr) -V_{0}V_{0}^{T}\bigr\rrVert
_{F}^{2}\leq 2M^{\prime}{}^{2}
\varepsilon^{2},
\]
for sufficiently large $(n,p,s,r)$.
\end{corollary}

The result follows from the fact that the $2(p+r)\delta$ part of
Proposition~\ref{prop:pointestimate} is exponentially small; hence, it is dominated
by $%
M'^2\varepsilon^2$.

\section{Discussion}

\label{sec:disc}

In Section~\ref{sec:spec}, we state a result on posterior contraction
rate under the spectral norm.
A computationally efficient algorithm is developed in Section~\ref{sec:comp}
for the rank-one case. In Section~\ref{sec:furtherremark}, we discuss
the possibility of using a simpler prior for sparse PCA.

\subsection{Posterior convergence under spectral norm}

\label{sec:spec}

In proving Theorem~\ref{thmm:postsub}, there are some by-products
serving as
intermediate steps. The following theorem says that the posterior
distribution concentrates on the true covariance matrix under the spectral
norm, and the subspace projection matrix concentrates on the true subspace
projection matrix under the spectral norm. In addition, the posterior
distribution consistently estimates the rank of the true subspace. The
theorem holds under a slightly weaker assumption without assuming
$r\leq ms$.

\begin{thmm}
\label{thmm:postspec} Consider the same prior $\Pi$ and rate
$\varepsilon
$ as
in Theorem~\ref{thmm:postsub}. Assume $\varepsilon\rightarrow0$, $r\leq
m\log
p $ and $n\leq p^{m}$ for some constant $m>0$. Then there exists
$M_{\gamma
,K,m}>0$, such that for any $M>M_{\gamma,K,m}$, we have
%
\begin{eqnarray}
\sup_{\Sigma\in\mathcal{G}(p,s,r)}P_{\Sigma}^{n}\Pi \bigl( \Vert
\Gamma -\Sigma\Vert >M\varepsilon|X^{n} \bigr) &\leq&\exp \bigl(
-C_{(\gamma
,K,m,M)}n\varepsilon^{2} \bigr), \label{eq:noneedforlowerK}
\\
\sup_{\Sigma\in\mathcal{G}(p,s,r)}P_{\Sigma}^{n}\Pi \bigl( \bigl\Vert
VV^{T}-V_{0}V_{0}^{T}\bigr\Vert >M
\varepsilon|X^{n} \bigr)& \leq&\exp \bigl( -C_{(\gamma,K,m,M)}n
\varepsilon^{2} \bigr) ,\nonumber
\\
\sup_{\Sigma\in\mathcal{G}(p,s,r)}P_{\Sigma}^{n}\Pi \bigl( \xi \neq
r|X^{n} \bigr) &\leq&\exp \bigl( -C_{(\gamma,K,m,M)}n\varepsilon
^{2} \bigr), \label{eq:rankconsistency}
\end{eqnarray}
for some constant $C_{(\gamma,K,m,M)}$ only depending on $(\gamma,K,m,M)$.
\end{thmm}

\begin{remark}
It is not practical to assume that $K$ is known in Theorems \ref
{thmm:postsub} and \ref{thmm:postspec}. To weaken the
assumption, we can replace the prior in (\ref{eq:priorcons}) by
sampling $U\sim\operatorname{Unif}[L_n^{-1},L_n]$, for some sequence $L_n$
slowly grows to infinity as $n\rightarrow\infty$. Then the conclusions
of the two theorems still hold without knowing $K$.
\end{remark}

\begin{remark}
The posterior rate of convergence (\ref{eq:noneedforlowerK}) for
estimating the whole covariance matrix under the spectral norm does not
require the assumption $\Vert \theta_l\Vert ^2>K^{-1}$ in the definition of
$\mathcal{G}(p,s,r)$. To remove this assumption, we need a slightly
different prior with (\ref{eq:priorcons}) modified by sampling $U\sim
\operatorname{Unif}[0,(2K)^{1/2}]$. However, such modification may not
lead to
rank adaptation (\ref{eq:rankconsistency}) due to lack of eigengap,
which is critical for establishing the result in Theorem~\ref{thmm:postsub}.
\end{remark}

\begin{remark}
Results (\ref{eq:noneedforlowerK}) and (\ref{eq:rankconsistency})
together imply posterior convergence of the whole covariance matrix
under the Frobenius norm. This is because when $\xi=r$, we have
$\Vert \Gamma-\Sigma\Vert _F=\Vert V\Lambda V^T-V_0\Lambda_0V_0^T\Vert _F\leq\sqrt
{2r}\Vert V\Lambda V^T-V_0\Lambda_0V_0^T\Vert =\sqrt{2r}\Vert \Gamma-\Sigma\Vert $.
Hence the convergence rate for the loss $\Vert \Gamma-\Sigma\Vert _F$ is
$\sqrt
{r}\varepsilon=\sqrt{\frac{r^2s\log p}{n}}$.
\end{remark}

\citeauthor{pati12} (\citeyear{pati12}) consider estimating the whole
covariance matrix under spectral norm in a sparse factor model. Under
their assumption $rs\gtrsim\log p$, they obtain a posterior convergence
rate of $\sqrt{\frac{r^3s\log p}{n}}\sqrt{\log n}$ under the loss
function $\Vert \Gamma-\Sigma\Vert $, compared with our rate $\sqrt{\frac
{rs\log p}{n}}$.

Though an improvement over the result of \citeauthor{pati12}
(\citeyear
{pati12}), whether $\sqrt{\frac{rs\log p}{n}}$ is the optimal rate of
convergence for the loss functions $\Vert \Gamma-\Sigma\Vert $ and
$\Vert VV^T-V_0V_0^T\Vert $ is still an open problem. To the best of our
knowledge, the only minimax result addressing these two loss functions
for sparse PCA problem is in \citeauthor {cai13} (\citeyear {cai13}).
However, they
consider a different sparsity class, defined as
\[
\mathcal{G}_1(p,s,r)=\lleft\{ %
\begin{array} {c}
\displaystyle \Sigma=\sum_{l=1}^{r}\theta_{l}
\theta_{l}^{T}+I\dvtx \biggl|\bigcup_{1\leq l\leq
r}S_{0l}\biggr|
\leq s,\theta_{l}\in\mathbb{R}^{p},
\\
\displaystyle \theta_{l}^{T}\theta_{k}=0\mbox{ for }k\neq l,
\Vert \theta_{l}\Vert ^{2}\in \bigl(K^{-1},K
\bigr)%
\end{array} %
 \rright\}.
\]
Under the current setting, the results of \citeauthor {cai13}
(\citeyear {cai13}) can
be written as
\begin{eqnarray*}
\inf_{\hat{\Sigma}}\sup_{\Sigma\in\mathcal{G}_1(p,s,r)}P_{\Sigma
}^n
\Vert \hat{\Sigma}-\Sigma\Vert ^2&\asymp&\frac{s\log p}{n}+
\frac{r}{n},
\\
\inf_{\hat{V}}\sup_{\Sigma\in\mathcal{G}_1(p,s,r)}P_{\Sigma
}^n
\bigl\Vert \hat {V}\hat{V}^T-VV^T\bigr\Vert ^2&\asymp&
\frac{s\log p}{n}.
\end{eqnarray*}
Observe the relation that
\[
\mathcal{G}_1(p,s,r)\subset\mathcal{G}(p,s,r)\subset
\mathcal{G}_1(p,rs,r).
\]
Hence when $r\leq O(\log p)$, the minimax rates for the class $\mathcal
{G}(p,s,r)$ under both loss functions lie between $\frac{s\log p}{n}$
and $\frac{rs\log p}{n}$. We claim that the posterior convergence rate
obtained in Theorem~\ref{thmm:postspec} is optimal when $r\leq O(1)$.
For a growing $r$, it at most misses a factor of $r$.

\subsection{A computational strategy of the rank-one case}
\label{sec:comp}

Bayesian procedures using sparse priors are usually harder to compute because
the sampling procedure needs to mix all possible subsets. %
\citeauthor{castillo12} (\citeyear{castillo12}) develop an efficient
algorithm for computing exact posterior mean in the setting of Bayesian
sparse vector
estimation. They explore the combinatorial nature of the posterior mean
formula and show that it is sufficient to compute the coefficients of
some $p $th order polynomials. In this section, we use their idea to develop an
algorithm for computing approximate posterior mean for the single spike
model. In this rank-one case, there is no need for the prior to adapt
to the
rank. We do not need the prior to put constraint on the $l^{2}$ norm of the
eigenvector as in (\ref{eq:priorcons}). Thus we use the following simple
prior on the single spiked covariance:
\begin{longlist}[(1)]
\item[(1)] sample a cardinality $q$ according to the distribution $\pi$
supported on $\{1,2,\ldots,p\}$;

\item[(2)] given $q$, sample a support $S\subset\{1,2,\ldots,p\}$ with
cardinality $%
|S|=q$ uniformly from all ${p\choose q}$ subsets;

\item[(3)] given $S$, sample $\eta_{S}\sim N(0,I_{|S|\times|S|})$, let
$%
\eta^T=(\eta_S^T,\eta_{S^c}^T)=(\eta_S^T,0^T)$ and the covariance
matrix is
$\Gamma=\eta\eta^T+I$.
\end{longlist}

We choose $\pi$ to be $\pi(q)\propto\exp ( -\kappa q\log
p
) $
for some constant $\kappa>0$. We let $\varepsilon^{2}=\frac{s\log
p}{n}$ be
the minimax rate when $r=1$. The posterior distribution induced by the above
prior has the following desired property:

\begin{thmm}
\label{thmm:rankone} Assume $\varepsilon\rightarrow0$ and $n\leq p^{m}$ for
some constant $m>0$. Then there exists $M_{\gamma,K,m}>0$, such that for
any $M>M_{\kappa,K,m}$, we have
\[
\sup_{\Sigma\in\mathcal{G}(p,s,1)}P_{\Sigma}^{n}\Pi \bigl( \min \bigl\{
\Vert \eta-\theta\Vert ,\Vert \eta+\theta\Vert \bigr\} >M\varepsilon|X^{n}
\bigr) \leq \exp \bigl( -C_{(\kappa,K,m,M)}n\varepsilon^{2} \bigr) ,
\]
for some constant $C_{(\kappa,K,m,M)}>0$ only depending on $(\kappa
,K,m,M)$%
.
\end{thmm}

Note that the loss function is the $l^{2}$ norm, which is stronger than the
loss function used in Theorem~\ref{thmm:postsub}. The theorem above is proved
in the supplementary material [\citeauthor {gaosupp} (\citeyear
{gaosupp})]. We use the
posterior mean $%
\mathbb{E}_{\Pi}(\eta|X^{n})$ to estimate the spike $\theta$.

We present a way for computing $\mathbb{E}_{\Pi}(\eta|X^{n})$. Under the
rank-one situation, representation (\ref{eq:latent}) can be written as
%
\begin{equation}
X_{ij}=W_{i}\theta_{j}+Z_{ij},\qquad i=1,
\ldots,n, j=1,\ldots,p, \label{eq:latentone}
\end{equation}
with $Z_{ij}$ and $W_{i}$ following i.i.d. $N(0,1)$ for all $i$ and
$j$. Representation (\ref{eq:latentone}) resembles the Gaussian
sequence model
considered in Castillo and van~der Vaart (\citeyear{castillo12}). Following
their idea, the $j$th coordinate of $\mathbb{E}_{\Pi}(\eta|X^{n})$
can be
written as
\[
\mathbb{E}_{\Pi}\bigl(\eta_{j}|X^{n}\bigr)=
\frac{\int\eta_{j}\int
\prod_{i=1}^{n}\prod_{j=1}^{p}\phi(X_{ij}-W_{i}\eta_{j})\bolds
{\phi}%
(W^{n})\,dW^{n}\,d\Pi(\eta)}{\int\int\prod_{i=1}^{n}\prod_{j=1}^{p}\phi
(X_{ij}-W_{i}\eta_{j})\bolds{\phi}(W^{n})\,dW^{n}\,d\Pi(\eta)},
\]
where $\bolds{\phi}(W^{n})\,dW^{n}=\prod_{i=1}^{n}\phi
(W_{i})\,dW_{1}\cdots dW_{n}$ and $\phi$ is the density function of
$N(0,1)$. By
Fubini's theorem, we have
\[
\mathbb{E}_{\Pi}\bigl(\eta_{j}|X^{n}\bigr)=
\frac{\int
N_{n,j}(W^{n})\bolds
{\phi
}(W^{n})\,dW^{n}}{\int D_{n}(W^{n})\bolds{\phi}(W^{n})\,dW^{n}},
\]
where for each $W^{n}$,
\begin{eqnarray*}
&&D_{n}\bigl(W^{n}\bigr) \\
&&\qquad=\int\prod
_{i=1}^{n}\prod_{j=1}^{p}
\phi (X_{ij}-W_{i}\eta _{j})\,d\Pi(\eta)
\\
&&\qquad=\sum_{q=1}^{p}\frac{\pi(q)}{{{p\choose q}}}\sum
_{|S|=q}\prod_{j\notin
S} \Biggl
\{ \prod_{i=1}^{n}\phi(X_{ij})
\Biggr\} \prod_{j\in S} \Biggl\{ \int \prod
_{i=1}^{n}\phi(X_{ij}-W_{i}
\eta_{j})\phi(\eta_{j})\,d\eta _{j} \Biggr\} ,
\end{eqnarray*}
by the definition of the prior. In the same way,
\begin{eqnarray*}
N_{n,j}\bigl(W^{n}\bigr) &=&\int\eta_{j}\prod
_{i=1}^{n}\prod
_{k=1}^{p}\phi (X_{ik}-W_{i}
\eta_{k})\,d\Pi(\eta)
\\
&=&\sum_{q=1}^{p}\frac{\pi(q)}{{{p\choose q}}}\sum
_{|S|=q}\prod_{k\notin
S} \Biggl
\{ \prod_{i=1}^{n}\phi(X_{ik})
\Biggr\} \\
&&{}\times\prod_{k\in S,k\neq
j} \Biggl\{ \int\prod
_{i=1}^{n}\phi(X_{ik}-W_{i}
\eta_{k})\phi(\eta_{k})\,d\eta _{k} \Biggr\}
\\
&&{}\times\mathbb{I}\{j\in S\}\int\eta_{j}\prod
_{i=1}^{n}\phi (X_{ij}-W_{i}
\eta_{j})\phi(\eta_{j})\,d\eta_{j}.
\end{eqnarray*}
Define
\begin{eqnarray*}
f(X_{\cdot j}) &=&\prod_{i=1}^{n}
\phi(X_{ij}),
\\
h\bigl(X_{\cdot j},W^{n}\bigr) &=&\int\prod
_{i=1}^{n}\phi(X_{ij}-W_{i}
\eta _{j})\phi (\eta_{j})\,d\eta_{j},
\\
\xi\bigl(X_{\cdot j},W^{n}\bigr) &=&\int\eta_{j}
\prod_{i=1}^{n}\phi (X_{ij}-W_{i}
\eta_{j})\phi(\eta_{j})\,d\eta_{j}.
\end{eqnarray*}
Then we may rewrite $D_{n}(W^{n})$ and $N_{n,j}(W^{n})$ as
\[
D_{n}\bigl(W^{n}\bigr)=\sum_{q=1}^{p}
\frac{\pi(q)}{{{p\choose
q}}}C\bigl(q,W^{n}\bigr),\qquad  N_{n,j}
\bigl(W^{n}\bigr)=\sum_{q=1}^{p}
\frac{\pi(q)}{{{p\choose q}}}C_{j}\bigl(q,W^{n}\bigr).
\]
The critical fact observed by \citeauthor{castillo12} (\citeyear
{castillo12}%
) is that $C(q,W^{n})$ is the coefficient of $Z^{q}$ of the polynomial
\[
Z\mapsto\prod_{j=1}^{p} \bigl(
f(X_{\cdot j})+h\bigl(X_{\cdot
j},W^{n}\bigr)Z \bigr) ,
\]
and $C_{j}(q,W^{n})$ is the coefficient of $Z^{q}$ of the polynomial
\[
Z\mapsto\xi\bigl(X_{\cdot j},W^{n}\bigr)Z\prod
_{k\in\{1,\ldots,p\}\setminus
\{j\}} \bigl( f(X_{\cdot k})+h\bigl(X_{\cdot k},W^{n}
\bigr)Z \bigr).
\]
For a given $W^{n}$, the coefficients $\{C(q,W^{n})\}_{q}$ and $%
\{C_{j}(q,W^{n})\}_{(j,q)}$ can be computed efficiently. In the Gaussian
sequence model, there is no randomness by $W^{n}$, and the posterior mean
can be computed exactly by finding the coefficients of the above
polynomials. In the PCA case, we propose an approximation by first
drawing $%
W_{1}^{n},W_{2}^{n},\ldots,W_{T}^{n}$ i.i.d. from $N(0,I_{n\times n})$ and then
computing
%
\begin{eqnarray}\label{eq:montecarlo}
\hat{\theta}_{j}={\frac{1}{T}\sum_{t=1}^{T} \Biggl( \sum_{q=1}^{p}\frac{%
\pi(q)}{{{p\choose q}}}C\bigl(q,W_{t}^{n}\bigr) \Biggr) }\bigg/{\Biggl(\frac{1}{T}%
\sum_{t=1}^{T} \Biggl( \sum_{q=1}^{p}\frac{\pi(q)}{{{p\choose q}}}%
C_{j}\bigl(q,W_{t}^{n}\bigr) \Biggr) \Biggr)}
\nonumber
\\[-8pt]
\\[-8pt]
\eqntext{\mbox{for }j=1,2,\ldots,p.}
\end{eqnarray}
One set of coefficients takes at most $O(p^{2})$ steps to compute. Thus the
total computational complexity is $O(Tp^{3}+Tnp)$ for computing coefficients
of $O(Tp)$ polynomials and computing all the values of $f(X_{\cdot
j})$, $%
h(X_{.j},W^{n})$ and $\xi(X_{\cdot j},W^{n})$.

The above strategy can be directly generalized to the multiple rank
case. However, it only works for the following prior without the
ability for rank adaptation. To be specific, we assume the rank $r$ is
known. Then, the third step of the prior is modified as follows:
\begin{longlist}[(3)]
\item[(3)] Given $S$, sample an $|S|\times r$ matrix $A_{S}$, with each
entry i.i.d. $N(0,1)$. Let the matrix $A$ be defined as
\[
A= %
\pmatrix{ A_S
\cr
0 }.
\]
The covariance matrix is $\Gamma=AA^T+I$.
\end{longlist}

Note that instead of sampling an individual support $S_l$ for each
column of $A$, we sample a common support $S$ for all columns. When
$r\leq O(1)$, this will not be a problem because of the simple
observation $rs\asymp s$. The theoretical justification of the prior is
stated in Theorem~\ref{thmm:rankmultiple}. Denote the $j$th row of $A$
by $A_j^T$. Then the posterior mean has formula
\[
\mathbb{E}_{\Pi}\bigl(A_j|X^n\bigr)=
\frac{\int N_{n,j}(W^n)\phi
(W^n)\,dW^n}{\int
D_n(W^n)\phi(W^n)\,dW^n},
\]
where for each $W^n$, we have
\begin{eqnarray*}
D_n\bigl(W^n\bigr)&=&\sum_{q=1}^p
\frac{\pi(q)}{{p\choose q}}\sum_{|S|=q}\prod
_{j\notin S} \Biggl\{\prod_{i=1}^n
\phi(X_{ij}) \Biggr\}\\
&&{}\times \prod_{j\in
S} \Biggl\{
\int\prod_{i=1}^n\phi\bigl(X_{ij}-W_i^TA_j
\bigr)\phi(A_j)\,dA_j \Biggr\},
\end{eqnarray*}
and a similar formula for $N_{n,j}(W^n)$.
Note that the only difference from the rank-one case is the inner
product $W_i^TA_j$. The notation $W^n$ stands for $(W_1,\ldots,W_n)$,
where each $W_i$ is an $r$-dimensional standard Gaussian vector. A
similar formula holds for $N_{n,j}(W^n)$. Thus we can apply the same
Monte Carlo approximation (\ref{eq:montecarlo}) for $\mathbb{E}_{\Pi
}(A_j|X^n)$ as is done in the rank-one case.

In addition to our method, there are other methods proposed in the
literature. A Gaussian shrinkage prior for Bayesian PCA have been
developed by \citeauthor {bishop99a} (\citeyear {bishop99a,bishop99b}) in the classical setting, but it is not appropriate for
sparse PCA. More general shrinkage priors have been discussed in
\citeauthor
{polson10} (\citeyear {polson10}) and \citeauthor {bhattacharya12}
(\citeyear
{bhattacharya12}) for high-dimensional mean vector estimation. One can
extend the framework to sparse PCA and develop Gibbs sampling by taking
advantage of the latent representation (\ref{eq:latent}). We refer to
\citeauthor {pati12} (\citeyear {pati12}) and van~der Pas, Kleijn and van~der Vaart
(\citeyear {pas14}) for some
theoretical justifications of shrinkage priors.

\subsection{Further remarks on the prior}
\label{sec:furtherremark}

The prior we proposed in Section~\ref{sec:main} on the random
covariance matrix $\Gamma=AA^T+I$ imposes orthogonality on the columns
of $A$. The orthogonality constraint is convenient for creating an
eigengap between the spikes and the noise. This leads to the rank
adaptation (\ref{eq:rankconsistency}). One may wonder if a simpler
prior such as the one proposed in Section~\ref{sec:comp} without
orthogonality constraint would also lead to a desired eigengap.

The answer is negative in the current proof technique. Let us consider
the simplest case where the supports $S_{01},S_{02},\ldots,S_{0r}$ are
known and $S_{01}=S_{02}=\cdots=S_{0r}=S_0$. When the rank $r$ is not
known, it is necessary to sample $\xi$ according to some prior
distribution. Then, after sampling the rank $\xi$, we only need to
sample a $|S_0|\times\xi$ submatrix of $A$, with rows in $S_0$. Let us
denote the submatrix by $A_{S_0}$. Consider the prior distribution of
$A_{S_0}$ where each element follows i.i.d. $N(0,1)$. Assume $r\leq s$
so that we can also restrict $\xi<s$. It is easy to see that the $\xi
$th eigenvalue of the matrix $\Gamma=AA^T+I$ is $\lambda_{\min
}(A_{S_0}A_{S_0}^T)+1$. Hence the eigengap is $\lambda_{\min
}(A_{S_0}A_{S_0}^T)$. For rank adaptation (\ref{eq:rankconsistency}),
we need a positive eigengap $\lambda_{\min}(A_{S_0}A_{S_0}^T)>0$. By
nonasymptotic random matrix theory [\citeauthor {vershynin10}
(\citeyear {vershynin10})],
%
\begin{equation}
\Pi \bigl(\lambda_{\min}\bigl(A_{S_0}A_{S_0}^T
\bigr)>\sqrt{s}-\sqrt{\xi }-t |\xi \bigr)\geq1-2e^{-t^2/2}, \label{eq:VRMT}
\end{equation}
for any $t>0$. For $\sqrt{s}-\sqrt{\xi}-t>0$, $t$ cannot be larger than
$\sqrt{s}$, leading to a tail not smaller than $2\exp(-s/2)$. In order
that there is an eigengap under the posterior distribution, the desired
tail needed in the classical Bayes nonparametric theory [see
\citeauthor
{barron98} (\citeyear {barron98}) and \citeauthor {castillo08}
(\citeyear {castillo08})] is
$\exp(-Cn\varepsilon^2)=\exp(-Crs\log p)$ for some $C>0$. Hence the random
matrix theory tail in (\ref{eq:VRMT}) is not enough for our purpose,
and the current proof technique does not lead to the desired posterior
convergence for this simpler prior. One may consider a larger support
$S$ with $|S|\asymp rs\log p$ in the prior distribution, such that the
tail probability in (\ref{eq:VRMT}) is $\exp(-Crs\log p)$ for some
$C>0$. However, it can be shown that the prior does not have sufficient
mass around the truth.

Nonetheless, if we assume the rank is known and $r\leq O(1)$, then rank
adaptation is not needed. In this case, the prior in Section~\ref{sec:comp} leads to the desired posterior rate of convergence. Remember
$\varepsilon^2=\frac{s\log p}{n}$.

\begin{thmm}
\label{thmm:rankmultiple} Assume $\varepsilon\rightarrow0$, $n\leq
p^{m}$ and $r\leq m$ for
some constant $m>0$. Then there exists $M_{\gamma,K,m}>0$, such that for
any $M>M_{\kappa,K,m}$, we have
\[
\sup_{\Sigma\in\mathcal{G}(p,s,r)}P_{\Sigma}^{n}\Pi \bigl( \bigl
\llVert VV^T-V_0V_0^T\bigr\rrVert
_F >M\varepsilon|X^{n} \bigr) \leq \exp \bigl(
-C_{(\kappa,K,m,M)}n\varepsilon^{2} \bigr) ,
\]
for some constant $C_{(\kappa,K,m,M)}>0$ only depending on $(\kappa
,K,m,M)$%
.
\end{thmm}

It would be an interesting problem to consider whether new techniques
can be developed to prove optimal posterior rate of convergence for a
simpler prior when the rank $r$ is not known.

\section{Proofs}
\label{sec:proof}

The results of Theorems \ref{thmm:postsub} and \ref
{thmm:postspec} are
special cases for bounding
%
\begin{equation}
P_{\Sigma}^n\Pi\bigl(B|X^n\bigr)=P_{\Sigma}^n
\frac{N_n(B)}{D_n}, \label{eq:tobound}
\end{equation}
where $D_n=\int\frac{dP_{\Gamma}^n}{dP_{\Sigma}^n}(X^n)\,d\Pi(\Gamma)$
and $%
N_n(B)=\int_B\frac{dP_{\Gamma}^n}{dP_{\Sigma}^n}(X^n)\,d\Pi(\Gamma
)$ for
different $B$. To bound (\ref{eq:tobound}), it is sufficient to upper bound
the numerator $N_n(B)$ and lower bound the denominator $D_n$. Following
Barron, Schervish and Wasserman (\citeyear{barron99}) and \citeauthor
{ghosal00} (%
\citeyear{ghosal00}), this involves three steps:
\begin{longlist}[(1)]
\item[(1)] Show the prior $\Pi$ puts sufficient mass near the truth;
that is, we
need
\[
\Pi(K_{n})\geq\exp \bigl( -Cn\varepsilon^{2} \bigr) ,
\]
where $K_{n}= \{ \Gamma\dvtx \frac{\Vert \Gamma-\Sigma\Vert _{F}}{\lambda
_{\min
}(\Gamma)}\leq\varepsilon \} $.

\item[(2)] Choose an appropriate subset $\mathcal{F}$, and show the
prior is
essentially supported on $\mathcal{F}$ in the sense that
\[
\Pi\bigl(\mathcal{F}^{c}\bigr)\leq\exp \bigl( -Cn
\varepsilon^{2} \bigr).
\]
This controls the complexity of the prior. Note that it is sufficient to
have $\Pi(\mathcal{F}^{c}|X^{n})\leq\exp ( -Cn\varepsilon
^{2}
) $.

\item[(3)] Construct a testing function $\phi$ for the following
testing problem:
\[
H_{0}\dvtx \Gamma=\Sigma,\qquad H_{1}\dvtx \Gamma\in B\cap
\mathcal{F}.
\]
We need to control the testing error in the sense that
\[
P_{\Sigma}^{n}\phi\vee\sup_{\Gamma\in B\cap\mathcal{F}}P_{\Gamma
}^{n}(1-
\phi)\leq\exp \bigl( -Cn\varepsilon^{2} \bigr).
\]
\end{longlist}

Notice the constants $C$'s are different in the above three steps, and
should satisfy some constraints in the proof. Step 1 lower bounds the prior
concentration near the truth, which leads to a lower bound for $D_{n}$. In
its original form [\citeauthor{schwartz65} (\citeyear{schwartz65})], $K_{n}$
is taken to be a fixed neighborhood of the truth defined through
Kullback--Leibler divergence. Step 2 and
step 3 are mainly for upper bounding $N_{n}(B)$. The testing idea in
step 3
is initialized by \citeauthor{lecam73} (\citeyear{lecam73}) and %
\citeauthor{schwartz65} (\citeyear{schwartz65}). Step 2 goes back to %
\citeauthor{barron88} (\citeyear{barron88}), who proposes the idea to choose
an appropriate $\mathcal{F}$ to regularize the alternative hypothesis
in the
test; otherwise the testing function for step 3 may never exist; see
\citeauthor{lecam73} (\citeyear{lecam73}) and \citeauthor{barron89}
(\citeyear{barron89}).

We list key technical lemmas needed in the proof for all three steps as
follows. From now on, all capital letters $C$ with or without
subscripts are
absolute constants. They do not depend on other quantities unless otherwise
mentioned.

\begin{lemma}
\label{lem:denominator} Assume $\varepsilon\rightarrow0$. Then for any $b>0$,
we have
\[
P_{\Sigma}^{n} \bigl( D_{n}\leq\Pi(K_{n})
\exp \bigl( -(b+1)n\varepsilon ^{2} \bigr) \bigr) \leq\exp \bigl(
-4C_{2}b^{2}K^{-1}n\varepsilon ^{2}
\bigr) ,
\]
where $C_{2}>0$ is an absolute constant.
\end{lemma}

\begin{lemma}
\label{lem:priorcon} Assume $\varepsilon\rightarrow0$ and $r\vee\log
n\leq
m\log p$ for some $m>0$. Then we have
\[
\Pi(K_{n})\geq\exp \bigl( - ( \gamma+2+mC_{1}\log
K+mC_{1} ) n\varepsilon^{2} \bigr) ,
\]
with some absolute constant $C_{1}>0$.
\end{lemma}

Lemma~\ref{lem:denominator} lower bounds the denominator $D_{n}$. It
is a
general result for all Gaussian covariance matrix estimation problems.
Lemma~\ref{lem:priorcon} lower bounds $\Pi(K_{n})$ in step 1.

\begin{lemma}
\label{cor:sparsity} Let $S=S_{1}\cup\cdots\cup S_{\xi}$. Assume
$\varepsilon
\rightarrow0$. When $r\vee\log n\leq m\log p$ for some $m>0$, we have
\[
P_{\Sigma}^{n}\Pi \bigl( |S|>\operatorname{Ars}|X^{n} \bigr) \leq
\exp \biggl( -\frac{%
\gamma A}{8}n\varepsilon^{2} \biggr) +\exp \bigl(
-4C_{2}K^{-1}n\varepsilon ^{2} \bigr) ,
\]
for any $A>8\gamma^{-1} ( \gamma+4+mC_{1}\log K+mC_{1} ) $.
\end{lemma}

Lemma~\ref{cor:sparsity} establishes the sparse property of the prior
$\Pi$. It corresponds to step~2, where $\mathcal{F}$ is the sparse
subset $%
\{\Gamma\dvtx |S|\leq \operatorname{Ars}\}$. Note that the parameter space we consider requires
$\max_{1\leq l\leq r}|S_{0l}|\leq s$. The sparsity constraint in
$\mathcal{F}
$ is much weaker, which means $\mathcal{F}$ is larger than the parameter
space we consider. Since we only need $\mathcal{F}$ to control the
regularity of the parameters in the alternative for hypothesis testing in
step 3, the oversized $\mathcal{F}$ here does not cause a problem. In
many Bayes nonparametric problems, the
parameter space can be negligible compared with the set $\mathcal{F}$.
\citeauthor{zhao00} (\citeyear{zhao00}) provides an example where the
parameter space receives no prior probability, while the set $\mathcal
{F}$ receives prior probability close to one; see also \citeauthor
{vaart08} (\citeyear
{vaart08}).

\begin{lemma}
\label{lem:test1} Assume $\varepsilon\rightarrow0$. There exists some
constant $M_{A,K,m}$ depending only on $(A,K,m)$, such that for any
$M>M_{A,K,m}$, we
have a testing function $\phi$ satisfying
\[
P_{\Sigma}^{n}\phi\leq3\exp \biggl( -\frac
{C_{3}M^{2}}{8K^{2}}n
\varepsilon ^{2} \biggr)
\]
and
\[
\sup_{\Gamma\in \{ \Gamma\dvtx \Vert \Gamma-\Sigma
\Vert >M\varepsilon
,|S|\leq \operatorname{Ars} \} }P_{\Gamma}^{n} ( 1-
\phi ) \leq\exp \biggl( -%
\frac{C_{3}M}{8}n\varepsilon^{2}
\biggr).
\]
\end{lemma}

The existence of a test and its error rates in step 3 are established in
Lemma~\ref{lem:test1}. These lemmas prove Theorem~\ref%
{thmm:postspec}.

In order to prove Theorem~\ref{thmm:postsub}, we need to establish a stronger
testing procedure. Since we have the conclusion of Theorem~\ref
{thmm:postspec}%
, it is sufficient to consider the subset $\{\Gamma\dvtx \Vert \Sigma-\Gamma
\Vert \leq
M\varepsilon\}$. More specifically, we are going to test $\Sigma
=V_{0}\Lambda
_{0}V_{0}^{T}+I$ against the following alternative:
\[
\mathcal{H}_{1}=\bigl\{ %
\Gamma=V
\Lambda V^{T}+I\dvtx \bigl\Vert VV^{T}-V_{0}V_{0}^{T}
\bigr\Vert _{F}>M^{\prime
}\varepsilon, \xi=r,|S|\leq \operatorname{Ars}%
 \bigr\}.
\]
Note that $S=S_{1}\cup\cdots\cup S_{\xi}$ is the joint support. The existence
of the test is established by the following lemma.

\begin{lemma}
\label{lem:test2} Assume $\varepsilon\rightarrow0$, $r\vee\log n\leq
m\log p$
and $r\leq ms$ for some absolute constant $m>0$. There exists some
constant $M_{A,K,m}^{\prime}$ only depending on $(A,K,m)$, and for any
$M^{\prime
}>M_{A,K,m}^{\prime}$, we have a testing function $\phi$ such that
\[
P_{\Sigma}^{n}\phi\leq3\exp \bigl( -\tfrac{1}{8}C_{5}
\delta _{K}^{\prime}%
\bar{M}^{2}n
\varepsilon^{2} \bigr)
\]
and
\[
\sup_{\Gamma\in\mathcal{H}_{1}}P_{\Gamma}^{n}(1-\phi
)\leq2\exp \bigl( -C_{5}\delta_{K}^{\prime}
\bar{M}^{2}n\varepsilon^{2} \bigr) ,
\]
where $\bar{M}=2^{-3/2}K^{-1}M^{\prime}$, $\delta_{K}'$ only depending
on $K$,
and $C_{5}$ is an absolute constant.
\end{lemma}

We are going to develop the proofs in several parts. In Section~\ref%
{sec:proofmain}, we establish the main results based on the key lemmas
above. All key lemmas are proved in the later sections. In Section~\ref
{sec:priorcon}, we prove Lemma~\ref%
{lem:priorcon}, which is for the prior concentration (step 1). In
Section~\ref{sec:sparsity}, we prove Lemma~\ref{cor:sparsity} by showing that the
prior puts most mass on a sparse set (step 2). Sections~\ref{sec:testspec}
and~\ref{sec:testfrob} are devoted in proving Lemmas \ref{lem:test1}
and \ref{lem:test2}, respectively (step 3). The proof of Lemma~\ref{lem:denominator}
is stated in supplementary material
[\citeauthor{gaosupp} (\citeyear {gaosupp})].

\subsection{Proofs of the main results}
\label{sec:proofmain}

In this section we prove Theorems \ref{thmm:postsub} and~\ref{thmm:postspec}.
Since the proof of Theorem~\ref{thmm:postsub} depends on the
conclusion of
Theorem~\ref{thmm:postspec}, we prove the latter one first.

\subsubsection{Proof of Theorem \texorpdfstring{\protect\ref{thmm:postspec}}{4.1}}

We decompose the posterior by
\[
\Pi \bigl( \Vert \Gamma-\Sigma\Vert >M\varepsilon|X^{n} \bigr) \leq
\Pi \bigl( \Vert \Gamma-\Sigma\Vert >M\varepsilon,|S|\leq \operatorname{Ars}|X^{n}
\bigr) +\Pi \bigl( |S|>\operatorname{Ars}|X^{n} \bigr) ,
\]
where $S=S_{1}\cup\cdots\cup S_{\xi}$. By Lemma~\ref{cor:sparsity}, we have
\[
P_{\Sigma}^{n}\Pi \bigl( |S|>\operatorname{Ars}|X^{n} \bigr) \leq
\exp \bigl( -\gamma An\varepsilon^{2}/8 \bigr) +\exp \bigl(
-4C_{2}K^{-1}n\varepsilon ^{2} \bigr) ,
\]
for any $A>8\gamma^{-1}(\gamma+4+mC_{1}\log K+mC_{1})$. From now on, we
fix $A$ to be $A=9\gamma^{-1}(\gamma+4+mC_{1}\log K+mC_{1})$. Then it is
sufficient to bound
\[
P_{\Sigma}^{n}\Pi \bigl( \Vert \Gamma-\Sigma\Vert >M
\varepsilon,|S|\leq \operatorname{Ars}|X^{n} \bigr).
\]
Let $\phi$ be the testing function in Lemma~\ref{lem:test1}, and we have
\begin{eqnarray*}
&&P_{\Sigma}^{n}\Pi \bigl( \Vert \Gamma-\Sigma\Vert >M
\varepsilon,|S|\leq \operatorname{Ars}|X^{n} \bigr)
\\
&&\qquad\leq P_{\Sigma}^{n}\Pi \bigl( \Vert \Gamma-\Sigma\Vert >M
\varepsilon ,|S|\leq \operatorname{Ars}|X^{n} \bigr) \bigl\{ D_{n}>
\Pi(K_{n})\exp\bigl(-2n\varepsilon ^{2}\bigr) \bigr\} (1-
\phi)
\\
&&\qquad\quad{}+P_{\Sigma}^{n}\phi+P_{\Sigma}^{n} \bigl(
D_{n}<\Pi(K_{n})\exp \bigl(-2n\varepsilon^{2}
\bigr) \bigr).
\end{eqnarray*}

There are three terms on the right-hand side above. By Lemma~\ref
{lem:test1}%
, $P_{\Sigma}^{n}\phi\leq3\exp ( -\frac
{C_{3}M^{2}}{8K^{2}}n\varepsilon
^{2} ) $ for sufficiently large $M$. By Lemma~\ref
{lem:denominator}, we
have $P_{\Sigma}^{n} ( D_{n}<\Pi(K_{n})\exp(-2n\varepsilon
^{2}) )
\leq\exp ( -4C_{2}K^{-1}n\varepsilon^{2} ) $. Now it remains to
bound the first term. Let $H_{1}= \{ \Gamma\dvtx \Vert \Gamma-\Sigma
\Vert >M\varepsilon,|S|\leq \operatorname{Ars} \} $. We have
\begin{eqnarray*}
&&P_{\Sigma}^{n}\Pi \bigl( \Vert \Gamma-\Sigma\Vert >M
\varepsilon,|S|\leq \operatorname{Ars}|X^{n} \bigr) \bigl\{ D_{n}>
\Pi(K_{n})\exp\bigl(-2n\varepsilon ^{2}\bigr) \bigr\} (1-
\phi)
\\
&&\qquad=P_{\Sigma}^{n} \biggl( \biggl({\int_{H_{1}}\frac{dP_{\Gamma}^{n}}{%
dP_{\Sigma}^{n}}\,d\Pi(\Gamma)}/{D_{n}}\biggr) \bigl\{
D_{n}>\Pi(K_{n})\exp \bigl(-2n\varepsilon^{2}
\bigr) \bigr\} (1-\phi) \biggr)
\\
&&\qquad\leq\frac{\exp(2n\varepsilon^{2})}{\Pi(K_{n})}P_{\Sigma}^{n}\int_{H_{1}}%
\frac{dP_{\Gamma}^{n}}{dP_{\Sigma}^{n}}(1-\phi)\,d\Pi(\Gamma)
\\
&&\qquad=\frac{\exp(2n\varepsilon^{2})}{\Pi(K_{n})}\int_{H_{1}}P_{\Gamma
}^{n}(1-
\phi)\,d\Pi(\Gamma)
\\
&&\qquad\leq\frac{\exp(2n\varepsilon^{2})}{\Pi(K_{n})}\sup_{\Gamma\in
H_{1}}P_{\Gamma}^{n}(1-
\phi),
\end{eqnarray*}
which is bounded by $\exp ( -\frac{C_{3}M}{16}n\varepsilon
^{2}
) $
because $\sup_{\Gamma\in H_{1}}P_{\Gamma}^{n}(1-\phi)$ is upper bounded
by Lemma~\ref{lem:test1}, and $\Pi(K_{n})$ is lower bounded by Lemma~\ref%
{lem:priorcon} for sufficiently large~$M$. By summing up the error
probability, we have
\[
P_{\Sigma}^{n}\Pi \bigl( \Vert \Gamma-\Sigma\Vert >M
\varepsilon|X^{n} \bigr) \leq \exp \bigl( -C_{(\gamma,K,m,M)}n
\varepsilon^{2} \bigr) ,
\]
for some constant $C_{(\gamma,K,m,M)}$ only depending on $(\gamma,K,m,M)$.

To obtain the rest of the results, it is sufficient to prove
%
\begin{equation}
\bigl\{ \Vert \Gamma-\Sigma\Vert \leq M\varepsilon \bigr\} \subset\{\xi=r\} \label{eq:implyrank}
\end{equation}
and
%
\begin{equation}
\bigl\{ \Vert \Gamma-\Sigma\Vert \leq M\varepsilon \bigr\} \subset \bigl\{ \bigl\Vert
VV^{T}-V_{0}V_{0}^{T}\bigr\Vert \leq KM
\varepsilon \bigr\}. \label{eq:implysubspace}
\end{equation}
Note that
\[
\Gamma=\sum_{l=1}^{\xi}\eta_{l}
\eta_{l}^{T}+I,
\]
and the eigenvalues of the covariance $\Gamma$ are $(\Vert \eta
_{1}\Vert ^{2}+1,\ldots,\Vert \eta_{\xi}\Vert ^{2}+1,1,\ldots,1)$, where the first $\xi$
eigenvalues are in the range $[(2K)^{-1}+1,(2K)+1]$ as specified by the
prior. Similarly, the eigenvalues of the covariance $\Sigma$ are
$(\Vert \theta
_{1}\Vert ^{2}+1,\ldots,\Vert \theta_{r}\Vert ^{2}+1,1,\ldots,1)$, and the first $r$
eigenvalues are in the range $[K^{-1}+1,K+1]$. Suppose $r<\xi$, let
$v\in\operatorname{span}(V)\cap\operatorname{span}(V_0)^{\perp}$
and $\Vert v\Vert =1$. Then
$v^T\Sigma v=1$ and $v^T\Gamma v\geq\lambda_{\xi}(\Gamma)\geq
1+(2K)^{-1}$, which contradicts $\Vert \Gamma-\Sigma\Vert \leq M\varepsilon$. The
same argument leads to contradiction when $r>\xi$. Thus we must have
$\xi=r$ when $\Vert \Gamma-\Sigma\Vert \leq M\varepsilon$.

Finally, (\ref{eq:implysubspace}) is an immediate consequence of the
Davis--Kahan sin-theta theorem (Lemma~\ref{lem:sintheta}). It is easy
to check that the eigengap $\delta$ in Lemma~\ref{lem:sintheta} is $K^{-1}$.

\subsubsection{Proof of Theorem \texorpdfstring{\protect\ref{thmm:postsub}}{3.1}}

With the results from Lemma~\ref{cor:sparsity} and Theorem~\ref
{thmm:postspec}%
, we decompose the posterior distribution as follows:
\begin{eqnarray*}
&&\Pi \bigl( \bigl\Vert VV^{T}-V_{0}V_{0}\bigr\Vert
_{F}>M^{\prime}\varepsilon |X^{n} \bigr)
\\
&&\qquad\leq\Pi \bigl( \bigl\Vert VV^{T}-V_{0}V_{0}\bigr\Vert
_{F}>M^{\prime}\varepsilon ,\Vert \Gamma -\Sigma\Vert \leq M
\varepsilon,|S|\leq \operatorname{Ars}|X^{n} \bigr)
\\
&&\qquad\quad{}+\Pi \bigl( \Vert \Gamma-\Sigma\Vert >M\varepsilon|X^{n} \bigr) +
\Pi \bigl( |S|>\operatorname{Ars}|X^{n} \bigr)
\\
&&\qquad\leq\Pi \bigl( \bigl\Vert VV^{T}-V_{0}V_{0}\bigr\Vert
_{F}>M^{\prime}\varepsilon, \xi =r,|S|\leq \operatorname{Ars}|X^{n}
\bigr)
\\
&&\qquad\quad{}+\Pi \bigl( \Vert \Gamma-\Sigma\Vert >M\varepsilon|X^{n} \bigr) +
\Pi \bigl( |S|>\operatorname{Ars}|X^{n} \bigr),
\end{eqnarray*}
where the last inequality is due to (\ref{eq:implyrank}).
Note that the later two terms converge to zero, as shown in Lemma~\ref
{cor:sparsity}
and Theorem~\ref{thmm:postspec}. Therefore, we only need to bound
\[
P_{\Sigma}^{n}\Pi \bigl(\bigl \Vert VV^{T}-V_{0}V_{0}
\bigr\Vert _{F}>M^{\prime
}\varepsilon , \xi=r,|S|\leq
\operatorname{Ars}|X^{n} \bigr).
\]
Remembering the definition of $\mathcal{H}_1$,
then, by Lemma~\ref{lem:test2}, there exists a testing function $\phi$
for $%
\mathcal{H}_{1}$ with the desired error bound. Using a similar argument
as in the
proof of Theorem~\ref{thmm:postspec}, we have established Theorem~\ref
{thmm:postsub}.

\subsection{The prior concentration of \texorpdfstring{$\Pi$}{Pi}}
\label{sec:priorcon}

We prove Lemma~\ref{lem:priorcon} in this
section. The main strategy for proving Lemma~\ref{lem:priorcon} is to
explore the structure of the prior. Specifically, since the prior $\Pi
$ is
defined by a sampling procedure for $\eta_{l+1}$ conditioning on
$\operatorname
{span%
}\{\eta_{1},\ldots,\eta_{l}\}$, we need to take advantage of this
feature by
using the chain rule and conditional independence.

\begin{pf*}{Proof of Lemma~\ref{lem:priorcon}} Since $\lambda_{\min
}(\Gamma
)\geq1$, we have
\[
\frac{\Vert \Gamma-\Sigma\Vert _{F}}{\lambda_{\min}(\Gamma)}\leq\Vert \Gamma -\Sigma\Vert _{F}.
\]
Write
\begin{eqnarray*}
&&\Pi \bigl( \Vert \Gamma-\Sigma\Vert _{F}\leq\varepsilon \bigr)
\\
&&\qquad\geq\Pi \bigl( \Vert \Gamma-\Sigma\Vert _{F}\leq\varepsilon
|%
(S_{1},\ldots,S_{[p^{\gamma/2}]})=(S_{01},
\ldots,S_{0r},\varnothing ,\ldots,\varnothing) \bigr)
\\
&&\qquad\quad{}\times\Pi \bigl((S_{1},\ldots,S_{[p^{\gamma
/2}]})=(S_{01},
\ldots,S_{0r},\varnothing,\ldots,\varnothing) \bigr).
\end{eqnarray*}
The second term in the above product is
\begin{eqnarray*}
&&\Pi \bigl((S_{1},\ldots,S_{[p^{\gamma
/2}]})=(S_{01},
\ldots,S_{0r},\varnothing,\ldots,\varnothing) \bigr)
\\
&&\qquad\geq \prod_{l=1}^{r}\Pi (
S_{l}=S_{0l} ) \prod_{l=r+1}^{[p^{\gamma/2}]}
\biggl(1-\frac{1}{p^{\gamma+1}} \biggr)^{p}
\\
&&\qquad\geq \biggl(1-\frac{1}{p^{\gamma+1}} \biggr)^{p^{1+\gamma/2}} \prod
_{l=1}^{r} \biggl( \frac{1}{p^{\gamma+1}} \biggr)
^{|S_{0l}|}
\\
&&\qquad\geq \exp \bigl(-2p^{-\gamma/2} \bigr)p^{-rs(\gamma+1)}
\\
&&\qquad\geq\exp \bigl( -(\gamma+2)rs\log p \bigr)
\end{eqnarray*}
because $p^{-\gamma/2}$ is at a smaller order of $rs\log p$.
Then we lower bound
\[
\Pi \bigl(\Vert \Gamma-\Sigma\Vert _{F}\leq\varepsilon
|(S_{1},\ldots,S_{[p^{\gamma/2}]})=(S_{01},
\ldots,S_{0r},\varnothing ,\ldots,\varnothing) \bigr).
\]
When $(S_{1},\ldots,S_{[p^{\gamma/2}]})=(S_{01},\ldots,S_{0r},\varnothing
,\ldots,\varnothing)$, we have
\begin{eqnarray*}
\Vert \Gamma-\Sigma\Vert _{F} &=&\Biggl\llVert \sum
_{l=1}^{r}\eta_l \eta_l
^{T}-\sum_{l=1}^{r}
\theta_l \theta_l^{T}\Biggr\rrVert
_{F} \leq\sum_{l=1}^{r}\bigl\Vert
\eta_l \eta_l^{T}-\theta_l
\theta_l^{T}\bigr\Vert _{F}
\\
&=&\sum_{l=1}^{r}\bigl\Vert
\eta_{l,S_{0l}}\eta_{l,S_{0l}}^{T}-\theta _{l,S_{0l}}
\theta_{l,S_{0l}}^{T}\bigr\Vert _{F}
\\
&\leq&\sum_{l=1}^{r}\Vert
\eta_{l,S_{0l}}-\theta_{l,S_{0l}}\Vert \bigl(\Vert \theta _{l,S_{0l}}
\Vert _{\infty}+\Vert \eta_{l,S_{0l}}\Vert _{\infty} \bigr)
\\
&\leq&(\sqrt{2}+1)K^{1/2}\sum_{l=1}^{r}
\Vert \eta_{l,S_{0l}}-\theta _{l,S_{0l}}\Vert.
\end{eqnarray*}
We use the notation $G$ to represent the probability $G_{(S_{1},\ldots,S_{r})}$
defined in Section~\ref{sec:main}. By conditional independence, we have
\begin{eqnarray*}
&&\Pi \bigl(\Vert \Gamma-\Sigma\Vert _{F}\leq\varepsilon
|(S_{1},\ldots,S_{[p^{\gamma/2}]})=(S_{01},
\ldots,S_{0r},\varnothing ,\ldots,\varnothing) \bigr)
\\
&&\qquad=G \Biggl(\Biggl\llVert \sum_{l=1}^{r}
\eta\eta^{T}-\sum_{l=1}^{r}\theta
\theta^{T}\Biggr\rrVert _{F}\leq\varepsilon \Biggr)
\\
&&\qquad\geq G \Biggl((\sqrt{2}+1)K^{1/2}\sum_{l=1}^{r}
\Vert \eta _{l,S_{0l}}-\theta _{l,S_{0l}}\Vert \leq\varepsilon \Biggr)
\\
&&\qquad\geq G \bigl((\sqrt{2}+1)K^{1/2}\Vert \eta_{l,S_{0l}}-\theta
_{l,S_{0l}}\Vert \leq \varepsilon_{l}, l=1,\ldots,r \bigr),
\end{eqnarray*}
where $\sum_{l=1}^{r}\varepsilon_{l}\leq\varepsilon$. In particular, we choose
\[
\varepsilon_{i}=c(r,\varepsilon) (3\sqrt{2}K)^{i},\qquad i=1,
\ldots,r,
\]
with $c(r,\varepsilon)=\frac{2}{3}\varepsilon(3\sqrt{2}K)^{-r}$. Then as long
as $K\geq1$, we have
\[
K\sum_{i=1}^{l}\varepsilon_{i}
\leq\frac{1}{2}\varepsilon_{l+1}
\]
and
\[
\sum
_{i=1}^{r}\varepsilon_{i}\leq
\varepsilon.
\]
Define $\mathcal{T}_{l}=\bigcap_{i=1}^{l}\mathcal{U}_{i}$ with
\[
\mathcal{U}_{i}= \bigl\{ (\sqrt{2}+1)K^{1/2}\Vert
\eta_{i,S_{0i}}-\theta _{i,S_{0i}}\Vert \leq\varepsilon_{i}
\bigr\} \qquad \mbox{for }i=1,\ldots,r.
\]
Using the chain rule, we have
\[
G(\mathcal{T}_{r})=G(\mathcal{U}_{1})\prod
_{l=1}^{r-1}G(\mathcal {T}_{l+1}|%
\mathcal{T}_{l}).
\]
For each $G(\mathcal{T}_{l}|\mathcal{T}_{l-1})$, we present a lower bound
and prove it in the supplementary material [\citeauthor {gaosupp}
(\citeyear {gaosupp})].

\begin{proposition}
\label{prop:chainrule} For each $l=1,2,\ldots,r-1$, we have
\begin{eqnarray*}
G(\mathcal{T}_{l+1}|\mathcal{T}_{l})&\geq&\frac{c(r,\varepsilon
)}{2(2+\sqrt
{2}%
)e^{K/2}}(3
\sqrt{2}K)^{l+1}\\
&&{}\times \exp \biggl(-s\log\frac{(4\sqrt{2}+1)K^{1/2}}{
c(r,\varepsilon)}-s\log (2\sqrt{s}/3
) \biggr).
\end{eqnarray*}
Moreover, $G(\mathcal{U}_1)$ can be lower bounded by the above formula
with $%
l=0$.
\end{proposition}

Using this result, we have
\begin{eqnarray*}
G(\mathcal{U}_{1})\prod_{l=1}^{r-1}G(
\mathcal{T}_{l+1}|\mathcal{T}_{l}) &\geq& \biggl(
\frac{c(r,\varepsilon)}{2(2+\sqrt{2})e^{K/2}} \biggr)^{r} (3 
\sqrt{2}K )^{r(r+1)/2}
\\
&&{}\times\exp \biggl(-rs\log\frac{(4\sqrt{2}+1)K^{1/2}}{c(r,\varepsilon
)}%
-C_{1}rs\log
s \biggr)
\\
&\geq&\exp \biggl(-C_{1}r^{2}s\log K-C_{1}rs
\log\frac{1}{\varepsilon}%
-C_{1}rs\log s \biggr),
\end{eqnarray*}
for some absolute constant $C_{1}>0$ when $\frac{K}{\log K}\leq rs$.
Therefore, we have
\begin{eqnarray*}
&&\Pi \biggl(\frac{\Vert \Gamma-\Sigma\Vert _{F}}{\lambda_{\min}(\Gamma
)}\leq \varepsilon \biggr)\\
&&\qquad\geq\exp \biggl(-(
\gamma+2)rs\log p-C_{1}r^{2}s\log K-C_{1}rs\log
\frac{1}{\varepsilon}-C_{1}rs\log s \biggr).
\end{eqnarray*}
Since
\[
\varepsilon^{2}=\frac{rs\log p}{n},
\]
we have
\[
\Pi \biggl(\frac{\Vert \Gamma-\Sigma\Vert _{F}}{\lambda_{\min}(\Gamma
)}\leq \varepsilon \biggr)\geq\exp \bigl(- (
\gamma+2+mC_{1}\log K+mC_{1} )%
n
\varepsilon^{2} \bigr),
\]
under the assumption $r\vee\log n\leq m\log p$ for some constant
$m>0$.
\end{pf*}

\subsection{The sparsity of \texorpdfstring{$\Pi$}{Pi}}
\label{sec:sparsity}

We prove Lemma~\ref{cor:sparsity} in this section. The result is
implied by
the prior sparsity stated in the following lemma.

\begin{lemma}
\label{lem:sparsity} For the sparsity prior specified above, we have
for any
$A>0$,
\[
\Pi \bigl(|S_1\cup\cdots\cup S_\xi| \geq \operatorname{Ars} \bigr)\leq\exp
\biggl(-\frac
{A\gamma}{4}%
rs\log p \biggr).
\]
\end{lemma}

\begin{pf*}{Proof of Lemma~\ref{lem:sparsity}} First, we have
\[
\Pi \bigl(|S_{1}\cup\cdots\cup S_{\xi}|>\operatorname{Ars} \bigr)\leq\Pi
\bigl(|S_{1}\cup \cdots\cup S_{[p^{\gamma/2}]}|>\operatorname{Ars} \bigr).
\]
Note that there is a slight abuse of notation above. The $\{S_l\}
_{l=1}^{\xi}$ on the left-hand side are from $\{S_l\}
_{l=1}^{[p^{\gamma
/2}]}$ on the right-hand side by excluding those $S_l$ with $\eta_l=0$.
Let $B=|S_{1}\cup\cdots\cup S_{[p^{\gamma/2}]}|$. Note
that $B$ is a Binomial random variable with parameter $\alpha$
satisfying $\alpha\leq p^{-1-\gamma/2}$.
Therefore,
\begin{eqnarray*}
\Pi (B>\operatorname{Ars} ) &\leq&\sum_{k=[\operatorname{Ars}]}^{p}{
\pmatrix{p
\cr
k}}\alpha ^{k}(1-\alpha)^{p-k} \leq\sum
_{k=[\operatorname{Ars}]}^{p}{\pmatrix{p
\cr
k}}\alpha ^{k}
\\
&\leq&\sum_{k=[\operatorname{Ars}]}^{p}\exp (k\log p )
\bigl(p^{-1-\gamma
/2} \bigr)^{k}
\\
&\leq&\sum_{k=[\operatorname{Ars}]}^{p}\exp \biggl(-k
\frac{\gamma}{2}\log p \biggr) \leq\exp \biggl(-\frac{A\gamma}{4}rs\log p
\biggr).
\end{eqnarray*}
Thus the proof is complete.
\end{pf*} %

Now we are ready to prove Lemma~\ref{cor:sparsity} by upper bounding the
numerator and lower bounding the denominator of $\Pi (|S_{1}\cup
\cdots\cup
S_{\xi}|>\operatorname{Ars}|X )$. This can be done by combining the results of
Lemmas %
\ref{lem:sparsity}, \ref{lem:denominator} and \ref{lem:priorcon}.

\begin{pf*}{Proof of Lemma~\ref{cor:sparsity}} Since $D_{n}=\int\frac{
dP_{\Gamma}^{n}}{dP_{\Sigma}^{n}}(X)\,d{\Pi}(\Gamma)$ and
$K_{n}=\break
\{
\frac{\Vert \Gamma-\Sigma\Vert _{F}}{\lambda_{\min}(\Gamma)}\leq\varepsilon
 \} $, we have
\begin{eqnarray*}
&&P_{\Sigma}^{n}\Pi\bigl (|S_{1}\cup\cdots\cup
S_{\xi}|>\operatorname{Ars}|X \bigr)
\\
&&\quad\leq P_{\Sigma}^{n}\Pi \bigl(|S_{1}\cup\cdots\cup
S_{\xi}|>\operatorname{Ars}|X \bigr) 
 \bigl\{ D_{n}\geq
\Pi(K_{n})\exp \bigl(-(b+1)n\varepsilon^{2} \bigr) \bigr\}
\\
&&\quad\quad{}+P_{\Sigma}^{n} \bigl\{ D_{n}\leq
\Pi(K_{n})\exp \bigl(-(b+1)n\varepsilon^{2}%
 \bigr) \bigr\}
\\
&&\qquad\leq\frac{\exp ((b+1)n\varepsilon^{2} )}{\Pi
(K_{n})}P_{\Sigma
}^{n}\int_{|S_{1}\cup\cdots\cup S_{\xi}|>\operatorname{Ars}}
\frac{dP_{\Gamma}^{n}}{%
dP_{\Sigma}^{n}}(X)\,d{\Pi}(\Gamma)
\\
&&\qquad\quad{}+\exp \bigl(-4C_{2}K^{-1}b^{2}n
\varepsilon^{2} \bigr)
\\
&&\qquad\leq\exp \bigl((b+1)n\varepsilon^{2} \bigr)\frac{\Pi(|S_{1}\cup
\cdots\cup
S_{\xi}|>\operatorname{Ars})}{\Pi(K_{n})}+\exp
\bigl(-4C_{2}K^{-1}b^{2}n\varepsilon
^{2} \bigr),
\end{eqnarray*}
where we have used Lemma~\ref{lem:denominator}. Using Lemmas \ref%
{lem:sparsity} and Lemma~\ref{lem:priorcon}, we have
\[
\frac{\Pi(|S_{1}\cup\cdots\cup S_{\xi}|>\operatorname{Ars})}{\Pi(K_{n})}\leq\exp \biggl(-%
 \biggl(\frac{A\gamma}{4}- (
\gamma+2+mC_{1}\log K+mC_{1} ) \biggr)%
n
\varepsilon^{2} \biggr).
\]
Hence by choosing $b=1$, we have
\begin{eqnarray*}
&& P_{\Sigma}^{n}\Pi \bigl(|S_{1}\cup\cdots\cup
S_{\xi}|>\operatorname{Ars}|X \bigr) \\
&&\qquad\leq \exp%
 \biggl(- \biggl(
\frac{A\gamma}{4}- (\gamma+4+mC_{1}\log K+mC_{1} )
\biggr)%
n\varepsilon^{2} \biggr)
\\
&&\qquad\quad{}+\exp \bigl(-4C_{2}K^{-1}n\varepsilon^{2}
\bigr).
\end{eqnarray*}
The conclusion then follows by letting $A>8\gamma^{-1} (\gamma
+4+mC_{1}\log K+mC_{1} )$.
\end{pf*} 

\subsection{Testing in spectral norm}
\label{sec:testspec}

We prove Lemma~\ref{lem:test1} in this section. Because of the
constraint $%
|S_{1}\cup\cdots\cup S_{\xi}|\leq \operatorname{Ars}$, we can break the testing problem into
many low-dimensional testing problems. Then a final test can be constructed
by combining the small tests. The following lemma establishes the existence
of such a low-dimensional test and bounds its error probability.

\begin{lemma}
\label{lem:spectest} For the random variable $Y^n=(Y_1,\ldots,Y_n)$ in
$\mathbb{%
R}^d$ and any $M>0$, there exists a testing function $\phi$, such that
\begin{eqnarray*}
&& P_{\bar{\Sigma}}^n\phi\bigl(Y^n\bigr)\leq\exp
\biggl(C_3\,d-\frac
{C_3M^2}{4\Vert \bar
{\Sigma}%
\Vert ^2}n\varepsilon^2 \biggr) + 2
\exp \bigl(C_3d-C_3M^{1/2}n \bigr),
\\
&&\sup_{\{\bar{\Gamma}\dvtx \Vert \bar{\Gamma}-\bar{\Sigma}\Vert >M\varepsilon\}
}P_{\bar
{\Gamma%
}}^n \bigl(1-\phi
\bigl(Y^n\bigr) \bigr)\\
&&\qquad\leq\exp \biggl(C_3d-
\frac{C_3Mn\varepsilon
^2}{4}%
\max \biggl\{1, \frac{M}{(M^{1/2}+2)^2\Vert \bar{\Sigma}\Vert ^2} \biggr\} \biggr),
\end{eqnarray*}
with some absolute constant $C_3>0$.
\end{lemma}

Notice $\bar{\Sigma}$ is a general $d\times d$ covariance matrix for
some $d$. It will be specified in the proof of Lemma~\ref{lem:test1}.
To prove Lemma~\ref{lem:spectest}, we need the following random matrix
inequality. Its
proof is given in the supplementary material [\citeauthor {gaosupp}
(\citeyear {gaosupp})].

\begin{lemma}
\label{lem:specconcen} Let $Y_1,\ldots,Y_n$ be i.i.d. from $N(0,\bar
{\Sigma})$,
where $\bar{\Sigma}$ is a $d\times d$ covariance matrix. Let $\hat
{\Sigma}=%
\frac{1}{n}\sum_{i=1}^nY_iY_i^T$ be the sample covariance matrix, and then
there is an absolute constant $C_3>0$, such that for any $t>0$,
\[
P_{\bar{\Sigma}}^n \bigl(\Vert \hat{\Sigma}-\bar{\Sigma}\Vert >t\Vert
\bar {\Sigma }\Vert \bigr)%
\leq\exp \bigl(-C_3 \bigl(-d+n\bigl(t
\wedge t^2\bigr) \bigr) \bigr).
\]
\end{lemma}

\begin{pf*}{Proof of Lemma~\ref{lem:spectest}}
Denote the alternative set
by $%
H_{1}=\{\bar{\Gamma}\dvtx \Vert \bar{\Gamma}-\bar{\Sigma}\Vert >M\varepsilon\}$,
and then
it will have following decomposition:
\[
H_{1}\subset\bigcup_{j=0}^{\infty}H_{1j},
\]
where
\[
H_{10}= \bigl\{ \Vert \bar{\Gamma}-\bar{\Sigma}\Vert >M\varepsilon,
\Vert \bar {\Gamma }\Vert \leq \bigl(M^{1/2}+2\bigr)\Vert \bar{\Sigma}
\Vert \bigr\} ,
\]
and for $j\geq1$,
\[
H_{1j}= \bigl\{ \bigl(M^{1/2}+2\bigr) \bigl(M
\varepsilon^{2} \bigr)^{-(j-1)/2}\Vert \bar {\Sigma}%
\Vert <\Vert \bar{\Gamma}\Vert \leq\bigl(M^{1/2}+2\bigr) \bigl(M
\varepsilon^{2} \bigr)^{-j/2}\Vert \bar {%
\Sigma}
\Vert \bigr\}.
\]
We divide the alternative set into pieces so that the spectral norm of
$\bar{%
\Gamma}$ is bounded in each piece. For the prior in Section~\ref{sec:main},
this is not needed because the prior only samples a random covariance matrix
with bounded spectrum. However, the prior in Section~\ref{sec:comp}
does not
impose a bounded spectrum constraint. The strategy for dividing the
alternative set is general for both cases.

We test each alternative hypothesis separately and then combine the
test and
use the union bound to control the error. To test against $H_{10}$, we use
\[
\phi_{0}=\mathbb{I} \Biggl\{ \Biggl\llVert \frac{1}{n}\sum
_{i=1}^{n}Y_{i}Y_{i}^{T}-
\bar{%
\Sigma}\Biggr\rrVert >M\varepsilon/2 \Biggr\}.
\]
To test against $H_{1j}$, we use
\[
\phi_{j}=\mathbb{I} \Biggl\{ \Biggl\Vert \frac{1}{n}\sum
_{i=1}^{n}Y_{i}Y_{i}^{T}
\Biggr\Vert>\frac{M^{1/2}+2}{2}\Vert \bar{\Sigma}\Vert \bigl(M
\varepsilon^{2} \bigr)%
^{-(j-1)/2} \Biggr\}.
\]
From Lemma~\ref{lem:specconcen}, we have
\[
P_{\bar{\Sigma}}^{n}\phi_{0}\leq\exp
\biggl(C_{3}d-\frac
{C_{3}M^{2}}{4\Vert %
\bar{\Sigma}\Vert ^{2}}n\varepsilon^{2} \biggr)
\]
and
\begin{eqnarray*}
P_{\bar{\Sigma}}^{n}\phi_{j} &\leq&P_{\bar{\Sigma}}^{n}
\Biggl\{ \Biggl\llVert \frac{1}{n}\sum_{i=1}^{n}Y_{i}Y_{i}^{T}-
\bar{\Sigma}\Biggr\rrVert +\Vert \bar{ 
\Sigma}\Vert >\frac{M^{1/2}+2}{2}
\Vert \bar{\Sigma}\Vert \bigl(M\varepsilon ^{2} \bigr)%
^{-(j-1)/2}
\Biggr\}
\\
&\leq&P_{\bar{\Sigma}}^{n} \Biggl\{ \Biggl\llVert
\frac{1}{n}%
\sum_{i=1}^{n}Y_{i}Y_{i}^{T}-
\bar{\Sigma}\Biggr\rrVert >\frac
{M^{1/2}}{2}\Vert %
\bar{\Sigma}\Vert
\bigl(M\varepsilon^{2} \bigr)^{-(j-1)/2} \Biggr\}
\\
&\leq&\exp \bigl(C_{3}d-C_{3}M^{1-j/2}n
\varepsilon^{-(j-1)} \bigr).
\end{eqnarray*}
Next, we control the type II error. For any $\bar{\Gamma}\in H_{10}$, we
have
\begin{eqnarray*}
P_{\bar{\Gamma}}^{n}(1-\phi_{0}) &\leq&P_{\bar{\Gamma}}^{n}
\Biggl\{ \Vert \bar{%
\Gamma}-\bar{\Sigma}\Vert -\Biggl\llVert
\frac{1}{n}\sum_{i=1}^{n}Y_{i}Y_{i}^{T}-%
\bar{\Gamma}\Biggr\rrVert <M\varepsilon/2 \Biggr\}
\\
&\leq&P_{\bar{\Gamma}}^{n} \Biggl\{ \Biggl\llVert
\frac{1}{n}%
\sum_{i=1}^{n}Y_{i}Y_{i}^{T}-
\bar{\Gamma}\Biggr\rrVert >M\varepsilon /2 \Biggr\}
\\
&\leq&P_{\bar{\Gamma}}^{n} \Biggl\{ \Biggl\llVert
\frac{1}{n}%
\sum_{i=1}^{n}Y_{i}Y_{i}^{T}-
\bar{\Gamma}\Biggr\rrVert >\Vert \bar {\Gamma }\Vert \frac{%
M\varepsilon}{2(M^{1/2}+2)\Vert \bar{\Sigma}\Vert } \Biggr\}
\\
&\leq&\exp \biggl(C_{3}d-\frac{C_{3}M^{2}}{4(M^{1/2}+2)^{2}\Vert \bar
{\Sigma
}%
\Vert ^{2}}n\varepsilon^{2}
\biggr).
\end{eqnarray*}
For any $H_{1j}$, we have
\begin{eqnarray*}
P_{\bar{\Gamma}}^{n}(1-\phi_{j}) &\leq&P_{\bar{\Gamma}}^{n}
\Biggl\{ \Vert \bar{%
\Gamma}\Vert -\Biggl\llVert \frac{1}{n}\sum
_{i=1}^{n}Y_{i}Y_{i}^{T}-
\bar {\Gamma }%
\Biggr\rrVert <\frac{M^{1/2}+2}{2}\Vert \bar{\Sigma}\Vert
\bigl(M\varepsilon ^{2} \bigr) 
^{-(j-1)/2} \Biggr\}
\\
&\leq&P_{\bar{\Gamma}}^{n} \Biggl\{ \Biggl\llVert
\frac{1}{n}%
\sum_{i=1}^{n}Y_{i}Y_{i}^{T}-
\bar{\Gamma}\Biggr\rrVert >\frac
{M^{1/2}+2}{2}\Vert %
\bar{\Sigma}\Vert
\bigl(M\varepsilon^{2} \bigr)^{-(j-1)/2} \Biggr\}
\\
&\leq&P_{\bar{\Gamma}}^{n} \Biggl\{ \Biggl\llVert
\frac{1}{n}%
\sum_{i=1}^{n}Y_{i}Y_{i}^{T}-
\bar{\Gamma}\Biggr\rrVert >\Vert \bar {\Gamma}%
\Vert M^{1/2}
\varepsilon/2 \Biggr\}
\\
&\leq&\exp \biggl(C_{3}d-\frac{C_{3}M}{4}n\varepsilon^{2}
\biggr).
\end{eqnarray*}
Now we combine the test by $\phi=\max_{0\leq j\leq\infty}\phi
_{j}$. The
error of the combined test can be bounded by
\begin{eqnarray*}
P_{\bar{\Sigma}}^{n}\phi&\leq&\sum_{j=0}^{\infty}P_{\Sigma
}^{n}
\phi_{j}
\\
&\leq&\exp \biggl(C_{3}d-\frac{C_{3}M^{2}}{4\Vert \bar{\Sigma
}\Vert ^{2}}n\varepsilon
^{2} \biggr)+\exp (C_{3}d )\sum
_{j=1}^{\infty}\exp \biggl(%
-C_{3}Mn
\varepsilon \biggl(\frac{1}{M^{1/2}\varepsilon} \biggr)^{j} \biggr)
\\
&\leq&\exp \biggl(C_{3}d-\frac{C_{3}M^{2}}{4\Vert \bar{\Sigma
}\Vert ^{2}}n\varepsilon
^{2} \biggr)+\exp (C_{3}d )\sum
_{j=1}^{\infty}\exp \biggl(%
-jC_{3}Mn
\varepsilon \biggl(\frac{1}{M^{1/2}\varepsilon} \biggr) \biggr)
\\
&\leq&\exp \biggl(C_{3}d-\frac{C_{3}M^{2}}{4\Vert \bar{\Sigma
}\Vert ^{2}}n\varepsilon
^{2} \biggr)+2\exp \bigl(C_{3}d-C_{3}M^{1/2}n
\bigr)
\end{eqnarray*}
and
\begin{eqnarray*}
P_{\bar{\Gamma}}^{n}(1-\phi) &\leq&P_{\bar{\Gamma}}^{n}
\min_{j}(1-\phi _{j})
\\
&\leq&\exp \biggl(C_{3}d-\frac{C_{3}Mn\varepsilon^{2}}{4}\max \biggl\{ 1,
\frac{M%
}{(M^{1/2}+2)^{2}\Vert \bar{\Sigma}\Vert ^{2}} \biggr\} \biggr).
\end{eqnarray*}
Thus the proof is complete.
\end{pf*} 

To prove Lemma~\ref{lem:test1}, we combine the small tests and control the
error by union bound.

\begin{pf*}{Proof of Lemma~\ref{lem:test1}} We denote the alternative
set by
\[
H_{1}= \bigl\{ \Gamma\dvtx \Vert \Gamma-\Sigma\Vert >M
\varepsilon,|S_{1}\cup \cdots\cup S_{\xi}|<\operatorname{Ars} \bigr\}.
\]
Define $S=S_{1}\cup\cdots\cup S_{\xi}$ and $S_{0}=S_{01}\cup\cdots\cup S_{0r}$.
We decompose $H_{1}$ by
\[
H_{1}\subset\bigcup_{B\dvtx |B|<\operatorname{Ars}}H_{1,B},
\]
where $H_{1,B}= \{ \Gamma\dvtx \Vert \Gamma-\Sigma\Vert >M\varepsilon
,S=B
\} $.
Define $\bar{B}=S\cup S_{0}$, and it is easy to see that
\[
\Vert \Gamma-\Sigma\Vert =\Vert \bar{\Gamma}-\bar{\Sigma}\Vert ,
\]
where
\[
\bar{\Gamma}=\sum_{l=1}^{\xi}
\eta_{l,\bar{B}}\eta_{l,\bar
{B}}^{T}+I,\qquad \bar{\Sigma}=\sum
_{l=1}^{r}\theta_{l,\bar{B}}
\theta_{l,\bar{B}}^{T}+I.
\]
Thus it is sufficient to test the following sub-problem in $\mathbb{R}^{
\bar{B}}$ for each $B$:
\[
H_{0}^{\prime}\dvtx \bar{\Gamma}=\bar{\Sigma},\qquad
H_{1,B}^{\prime}\dvtx \Vert \bar {%
\Gamma}-\bar{\Sigma}
\Vert >M\varepsilon.
\]
By Lemma~\ref{lem:spectest}, there exists $\phi_{B}$ depending on the
observations $(Y_{1},\ldots,Y_{n})=(X_{1,\bar{B}},\ldots,X_{n,\bar{B}})$, such
that
\begin{eqnarray*}
P_{\Sigma}^{n}\phi_{B} &\leq&\exp
\biggl(C_{3}(A+1)rs-\frac
{C_{3}M^{2}}{%
4K^{2}}n\varepsilon^{2}
\biggr)\\
&&{}+2\exp \bigl(C_{3}(A+1)rs-C_{3}M^{1/2}n
\bigr)
\\
&\leq&3\exp \biggl(-C_{3} \biggl(\frac{M^{2}}{4K^{2}}-(A+1) \biggr)n
\varepsilon ^{2}%
 \biggr),
\\
\sup_{\Gamma\in H_{1,B}}P_{\Gamma}^{n}(1-
\phi_{B}) &\leq&\exp \biggl( 
C_{3}(A+1)rs-
\frac{C_{3}Mn\varepsilon^{2}}{4}\max \biggl\{ 1,\frac{M}{%
(M^{1/2}+2)^{2}K^{2}} \biggr\} \biggr)
\\
&\leq&\exp \biggl(-C_{3} \biggl(\frac{M}{4}-(A+1) \biggr)n
\varepsilon ^{2} \biggr).
\end{eqnarray*}
Then we combine the tests by $\phi=\max_{B}\phi_{B}$. By the union bound,
we have
\begin{eqnarray*}
P_{\Sigma}^{n}\phi&\leq& \Biggl(\sum
_{q=1}^{[\operatorname{Ars}]}{\pmatrix{p
\cr
q}} \Biggr)3\exp%
 \biggl(-C_{3} \biggl(\frac{M^{2}}{4K^{2}}-(A+1) \biggr)n
\varepsilon^{2} \biggr)
\\
&\leq&3\operatorname{Ars}\exp \biggl(\operatorname{Ars}\log\frac{ep}{\operatorname{Ars}} \biggr)\exp
\biggl(-C_{3} \biggl( 
\frac{M^{2}}{4K^{2}}-(A+1) \biggr)n
\varepsilon^{2} \biggr)
\\
&\leq&3\exp (2\operatorname{Ars}\log p )\exp \biggl(-C_{3} \biggl(
\frac
{M^{2}}{4K^{2}}%
-(A+1) \biggr)n\varepsilon^{2} \biggr)
\\
&\leq&3\exp \biggl(- \biggl(\frac{C_{3}M^{2}}{4K^{2}}-C_{3}(A+1)-2A
\biggr)%
n\varepsilon^{2} \biggr)
\end{eqnarray*}
and
\[
\sup_{\Gamma\in H_{1}}P_{\Gamma}^{n}(1-\phi)\leq\exp
\biggl(-C_{3} \biggl(%
\frac{M}{4}-(A+1) \biggr)n
\varepsilon^{2} \biggr).
\]
Hence the proof is complete by choosing sufficiently large $M$.
\end{pf*} %

\subsection{Testing in subspace distance \texorpdfstring{$d(\cdot,\cdot)$}{d(.,.)}}
\label{sec:testfrob}

We prove Lemma~\ref{lem:test2} in this section. At first thought, there
seems to be no obvious test for testing the subspace projection matrix
under the
distance $d(\cdot,\cdot)$ due to the complicated sparse and low-rank
structure. Our strategy is to break the alternative set into many
levels and
pieces. The goal is that for each piece, it is a low-dimensional small
testing problem in the following form:
\[
H_{0}\dvtx \bar{\Gamma}=\bar{\Sigma}, \qquad H_{1}\dvtx\bigl \Vert
\bar{\Gamma}-\bar {\Gamma}%
^{\prime}\bigr\Vert _{F}\leq
\delta_{K}\bigl\Vert \bar{\Sigma}-\bar{\Gamma}^{\prime}\bigr\Vert
_{F}.
\]
The small testing problem can be solved by considering the likelihood ratio
test. The error bound is stated in the following lemma. Its proof is given
in the supplementary material [\citeauthor {gaosupp} (\citeyear {gaosupp})].

\begin{lemma}
\label{lem:testfrob} Consider observations $Y^n=(Y_1,\ldots,Y_n)$ in
$\mathbb{R}%
^d$. There exist constants $\delta_K$ and $\delta^{\prime}_K$ only depending
on $K$, and a testing function $\phi$ such that
\begin{eqnarray*}
P_{\bar{\Sigma}}^n\phi\bigl(Y^n\bigr)&\leq&2\exp
\bigl(-C_5\delta_K^{\prime
}n\bigl\Vert \bar {
\Sigma%
}-\bar{\Gamma}^{\prime}\bigr\Vert _F^2
\bigr),
\\
\sup_{ \{\bar{\Gamma}\dvtx \Vert \bar{\Gamma}-\bar{\Gamma}^{\prime
}\Vert _F\leq
\delta_K\Vert \bar{\Sigma}-\bar{\Gamma}^{\prime}\Vert _F \}}P_{\bar
{\Gamma
}}^n%
 \bigl(1-
\phi\bigl(Y^n\bigr) \bigr)&\leq&2\exp \bigl(-C_5
\delta_K^{\prime}n\bigl\Vert \bar {\Sigma }-%
\bar{
\Gamma}^{\prime}\bigr\Vert _F^2 \bigr),
\end{eqnarray*}
where $C_5>0$ is an absolute constant.
\end{lemma}

We need a lemma to bound the covering number under different subspace
distances. We use $N(\delta,\mathcal{H},\rho)$ to denote the $\delta$
-covering number of $\mathcal{H}$ under the distance $\rho$. The
proof of
Lemma~\ref{lem:dcover} is given in the supplementary material
[\citeauthor
{gaosupp} (\citeyear {gaosupp})].

\begin{lemma}
\label{lem:dcover} For any $U\in\mathcal{U}(d,r)$, $R_1,R_2>0$ and
$\Lambda=%
\operatorname{diag}(\lambda_1,\ldots,\lambda_r)$ with $\lambda_1\geq
\lambda_2\geq
\cdots
\geq\lambda_r$, we have
\begin{eqnarray*}
&& \log N \bigl(R_1\varepsilon, \bigl\{V\in\mathcal{U}(d,r)\dvtx
d(U,V)\leq R_2\varepsilon \bigr\}, d_{\Lambda} \bigr)
\\
&&\qquad\leq dr\log \biggl(\frac{12\lambda_1 (R_2+1)}{R_1} \biggr) +r^2\log
\frac{6
\sqrt{r}}{\varepsilon}.
\end{eqnarray*}
\end{lemma}

Last but not least, we need the following sin-theta theorem to bound
the difference of subspaces by the difference of matrices.
%
\begin{lemma}[{[\citeauthor {davis70} (\citeyear {davis70})]}] \label
{lem:sintheta}
Consider symmetric matrices $F$ and~$\hat{F}$, with eigenvalue decomposition
\[
F=U_1D_1U_1^T+U_2D_2U_2^T,\qquad
\hat{F}=\hat{U}_1\hat{D}_1\hat {U}_1^T+
\hat{U}_2\hat{D}_2\hat{U}_2^T.
\]
If the eigenvalues $D_1$ are contained in an interval $(a,b)$, and the
eigenvalues $\hat{D}_2$ are excluded from the interval $(a-\delta
,b+\delta)$ for some $\delta>o$, then
\[
\bigl\llVert U_1U_1^T-
\hat{U}_1\hat{U}_1\bigr\rrVert _F\leq
\sqrt{2}\delta ^{-1}\llVert F-\hat{F}\rrVert _F
\]
and
\[
\bigl\llVert U_1U_1^T-
\hat{U}_1\hat{U}_1\bigr\rrVert \leq\delta^{-1}
\llVert F-\hat {F}\rrVert .
\]
\end{lemma}

\begin{pf*}{Proof of Lemma~\ref{lem:test2}}The proof has two major steps.

\textit{Step} 1: Decompose the alternative set into many levels and pieces.
We first decompose $\mathcal{H}_{1}$ by $\mathcal{H}_{1}\subset
\bigcup_{B\dvtx |F|\leq \operatorname{Ars}}H_{1,B}$,
where
\[
H_{1,B}= \bigl\{ \Gamma=V\Lambda V^{T}+I\dvtx \bigl\Vert
VV^{T}-V_{0}V_{0}^{T}\bigr\Vert
_{F}>M^{\prime}\varepsilon, \xi =r,S=B \bigr\}.
\]
Define $\bar{B}=B\cup S_{0}$ with $S_{0}=S_{01}\cup\cdots\cup S_{0r}$, and
\begin{eqnarray*}
V_{\bar{B}}&=& \bigl[ \Vert \eta_{1,\bar{B}}\Vert ^{-1}
\eta_{1,\bar
{B}},\ldots,\Vert \eta _{r,\bar{B}}\Vert ^{-1}
\eta_{r,\bar{B}} \bigr] , \\
 V_{0,\bar{B}}&=& \bigl[ \Vert
\theta_{1,\bar{B}}\Vert ^{-1}\theta_{1,\bar{B}},\ldots,\Vert
\theta_{r,\bar{B}
}\Vert ^{-1}\theta_{r,\bar{B}} \bigr].
\end{eqnarray*}
Note that both $V_{\bar{B}}$ and $V_{0,\bar{B}}$ are $|\bar
{B}|\times r$
matrices with $|\bar{B}|\leq(A+1)rs$, and $%
\Vert VV^{T}-V_{0}V_{0}^{T}\Vert _{F}=\Vert V_{\bar{B}}V_{\bar{B}}^{T}-V_{0,\bar
{B}}V_{0,%
\bar{B}}^{T}\Vert _{F}$. Then we can rewrite $H_{1,B}$ as
\[
H_{1,B}= \bigl\{ \Gamma=V\Lambda V^{T}+I\dvtx\bigl \Vert
V_{\bar{B}}V_{\bar
{B}}^{T}-V_{0,%
\bar{B}}V_{0,\bar{B}}^{T}
\bigr\Vert _{F}>M^{\prime}\varepsilon \bigr\} ,
\]
where we omit $\xi=r$ for simplicity of notation, and we consider both
$%
\Lambda$ and $\Lambda_{0}$ $r\times r$ diagonal matrices from now on.

Note that $\Vert \Lambda^{-1}\Vert _{\infty}\vee\Vert \Lambda\Vert _{\infty}\leq2K$
for any $\Gamma\in\operatorname{supp}(\Pi)$. We can show there
exists diagonal
matrices $\{\Lambda_{1},\ldots,\Lambda
_{T}\}\subset \{ \Lambda\dvtx \Vert \Lambda^{-1}\Vert _{\infty}\vee
\Vert \Lambda
\Vert _{\infty}\leq
2K \} $ such that
\[
\bigl\{ \Lambda\dvtx \bigl\Vert \Lambda^{-1}\bigr\Vert _{\infty}\vee\Vert
\Lambda\Vert _{\infty
}\leq 2K \bigr\} \subset\bigcup
_{t=1}^{T} \bigl\{ \Lambda\dvtx \Vert \Lambda-\Lambda
_{t}\Vert _{F}\leq \varepsilon \bigr\} ,
\]
where $\log T\leq r\log (12K\sqrt{r}\varepsilon^{-1} )$,
because we
regard $ \{
\Lambda\dvtx \Vert \Lambda^{-1}\Vert _{\infty}\vee\Vert \Lambda\Vert _{\infty}\leq2K
 \} $ as a
subset of $ \{ \Lambda\dvtx \Vert \Lambda\Vert _{F}\leq2K\sqrt{r} \}
$ so
that it is essentially a covering number calculation
in $\mathbb{R}^{r}$ as in \citeauthor{pollard90} (\citeyear{pollard90}).
We further decompose $H_{1,B}$ by $H_{1,B}\subset\bigcup_{t=1}^{T}H_{1,B,t}$%
, where
\[
H_{1,B,t}= \bigl\{ \Gamma=V\Lambda V^{T}+I\dvtx \bigl\Vert
V_{\bar{B}}V_{\bar
{B}}^{T}-V_{0,%
\bar{B}}V_{0,\bar{B}}^{T}
\bigr\Vert _{F}>M^{\prime}\varepsilon,\Vert \Lambda -\Lambda
_{t}\Vert _{F}\leq\varepsilon \bigr\} ,
\]
and decompose $H_{1,B,t}$ by $H_{1,B,t}\subset\bigcup_{j=1}^{\infty
}H_{1,B,t,j}$, where
\begin{eqnarray*}
&&H_{1,B,t,j}= \bigl\{ \Gamma=V\Lambda V^{T}+I\dvtx
jM^{\prime}\varepsilon <\bigl\Vert V_{\bar{%
B}}V_{\bar{B}}^{T}-V_{0,\bar{B}}V_{0,\bar{B}}^{T}
\bigr\Vert _{F}\leq (j+1)M^{\prime
}\varepsilon,\\
&&\hspace*{267pt}\Vert \Lambda-
\Lambda_{t}\Vert _{F}\leq\varepsilon \bigr\}.
\end{eqnarray*}

According to Lemma~\ref{lem:dcover}, there exists
\[
\{U_{1},\ldots,U_{N_{j}}\}\subset\mathcal{U}\bigl(|\bar{B}|,r\bigr)\cap
\bigl\{ U\dvtx jM^{\prime}\varepsilon<\bigl\Vert UU^{T}-V_{0,\bar{B}}V_{0,\bar
{B}}^{T}
\bigr\Vert _{F}\leq (j+1)M^{\prime}\varepsilon \bigr\} ,
\]
such that for some constant $\delta_{K}$ only depending on $K$,
\begin{eqnarray*}
&& \bigl\{ jM^{\prime}\varepsilon<\bigl\Vert V_{\bar{B}}V_{\bar
{B}}^{T}-V_{0,\bar
{B}%
}V_{0,\bar{B}}^{T}
\bigr\Vert _{F}\leq(j+1)M^{\prime}\varepsilon \bigr\}
\\
&&\qquad\subset\bigcup_{i=1}^{N_{j}} \bigl\{\bigl \Vert
V_{\bar{B}}\Lambda _{t}V_{\bar
{B}%
}^{T}-U_{i}
\Lambda_{t}U_{i}^{T}\bigr\Vert _{F}\leq (
\delta_{K}j\bar {M}-1 )%
\varepsilon \bigr\} ,
\end{eqnarray*}
where $\bar{M}=2^{-1/2}K^{-1}M^{\prime}$, and we may bound $N_{j}$ by
\begin{eqnarray*}
\log N_{j} &\leq&|\bar{B}|r\log \biggl(\frac{12\lambda_{1} (%
(j+1)M^{\prime}+1 )}{j\delta_{k}\bar{M}-1}
\biggr)+r^{2}\log \frac
{6\sqrt{%
r}}{\varepsilon}
\\
&\leq&(A+1)r^{2}s\log \bigl(48\sqrt{2}\delta_{K}^{-1}K
\bigr)+r^{2}\log (6%
\sqrt{r} )+\frac{1}{2}r^{2}
\log n,
\end{eqnarray*}
when we choose $M^{\prime}>\max \{ 2\sqrt{2}\delta
_{K}^{-1}K,\frac{1}{%
2} \} $. Using the triangle inequality, we have
\[
\bigl\Vert V_{\bar{B}}\Lambda V_{\bar{B}}^{T}-U_{i}
\Lambda _{t}U_{i}^{T}\bigr\Vert _{F}\leq
\bigl\Vert V_{\bar{B}}\Lambda_{t}V_{\bar{B}}^{T}-U_{i}
\Lambda _{t}U_{i}^{T}\bigr\Vert _{F}+\Vert
\Lambda-\Lambda_{t}\Vert _{F}.
\]
Therefore,
\begin{eqnarray*}
&& \bigl\{\bigl \Vert V_{\bar{B}}\Lambda_{t}V_{\bar{B}}^{T}-U_{i}
\Lambda _{t}U_{i}^{T}\bigr\Vert _{F}\leq (
\delta_{K}j\bar{M}-1 )\varepsilon ,\Vert \Lambda -
\Lambda_{t}\Vert _{F}\leq\varepsilon \bigr\}
\\
&& \qquad\subset\bigl\{ \bigl\Vert V_{\bar{B}}\Lambda V_{\bar{B}}^{T}-U_{i}
\Lambda _{t}U_{i}^{T}\bigr\Vert _{F}\leq (
\delta_{K}j\bar{M} )\varepsilon \bigr\}.
\end{eqnarray*}
By the sin-theta theorem (Lemma~\ref{lem:sintheta}), we have
\begin{eqnarray*}
\bigl\Vert U_{i}\Lambda_{t}U_{i}^{T}-V_{0,\bar{B}}
\Lambda_{0}V_{0,\bar{B}%
}^{T}\bigr\Vert _{F}
&\geq& 2^{-1/2}K^{-1}\bigl\Vert U_{i}U_{i}^{T}-V_{0,\bar{B}}V_{0,\bar{B}
}^{T}
\bigr\Vert \\
&\geq& 2^{-1/2}K^{-1}jM^{\prime}\varepsilon\geq j
\bar{M}\varepsilon.
\end{eqnarray*}
Hence
\begin{eqnarray*}
&& \bigl\{ \bigl\Vert V_{\bar{B}}\Lambda_{t}V_{\bar{B}}^{T}-U_{i}
\Lambda _{t}U_{i}^{T}\bigr\Vert _{F}\leq (
\delta_{K}j\bar{M}-1 )\varepsilon ,\Vert \Lambda -
\Lambda_{t}\Vert _{F}\leq\varepsilon \bigr\}
\\
&&\qquad\subset \bigl\{ \bigl\Vert V_{\bar{B}}\Lambda V_{\bar{B}}^{T}-U_{i}
\Lambda _{t}U_{i}^{T}\bigr\Vert _{F}\leq
\delta_{K}\bigl\Vert U_{i}\Lambda _{t}U_{i}^{T}-V_{0,\bar{B}%
}
\Lambda_{0}V_{0,\bar{B}}^{T}\bigr\Vert _{F} \bigr
\}.
\end{eqnarray*}
Our final decomposition is $H_{1,B,t,j}\subset
\bigcup_{i=1}^{N_{j}}H_{1,B,t,j,i}$, where
\begin{eqnarray*}
&& H_{1,B,t,j,i}= \bigl\{ \Gamma=V\Lambda V^{T}+I\dvtx\bigl \Vert
V_{\bar{B}}\Lambda V_{\bar{B%
}}^{T}-U_{i}
\Lambda_{t}U_{i}^{T}\bigr\Vert _{F}\\
&&\hspace*{98pt}\leq
\delta_{K}\bigl\Vert U_{i}\Lambda _{t}U_{i}^{T}-V_{0,\bar{B}}
\Lambda_{0}V_{0,\bar{B}}^{T}\bigr\Vert _{F} \bigr
\}.
\end{eqnarray*}

\textit{Step} 2: Combine tests from all levels and pieces. We have reduced
the original testing problem to the above small pieces for each $(B,t,j,i)$.
For each small piece, it is equivalent to the testing problem in Lemma~\ref%
{lem:testfrob}. Since we already know the coordinates $\bar{B}$, the
testing problem is on $\mathbb{R}^{\bar{B}}$. The observations in Lemma~\ref%
{lem:testfrob} is $(Y_{1},\ldots,Y_{n})=(X_{1,\bar{B}},\ldots,X_{n,\bar
{B}})$. The
triple $(\bar{\Sigma},\bar{\Gamma}^{\prime},\bar{\Gamma})$ in Lemma~\ref%
{lem:testfrob} corresponds to $ (V_{0,\bar{B}}\Lambda
_{0}V_{0,\bar
{B}%
}^{T}+I,U_{i}\Lambda_{t}U_{i}^{T}+I,V_{\bar{B}}\Lambda V_{\bar
{B}}^{T}+I%
 )$ for every $(B,t,j,i)$. Then by the conclusion of Lemma~\ref%
{lem:testfrob}, there exists a testing function $\phi_{B,t,j,i}$ with error
bounded by
\begin{eqnarray*}
P_{\Sigma}^{n}\phi_{B,t,j,i}&\leq& 2\exp
\bigl(-C_{5}\delta _{K}^{\prime
}n\bigl\Vert
U_{i}\Lambda_{t}U_{i}^{T}-V_{0,\bar{B}}
\Lambda_{0}V_{0,\bar{B}%
}^{T}\bigr\Vert _{F}^{2}
\bigr),
\\
\sup_{\Gamma\in H_{B,t,j,i}}P_{\Gamma}^{n}(1-
\phi_{B,t,j,i})&\leq& 2\exp%
 \bigl(-C_{5}
\delta_{K}^{\prime}n\bigl\Vert U_{i}\Lambda
_{t}U_{i}^{T}-V_{0,\bar
{B}%
}
\Lambda_{0}V_{0,\bar{B}}^{T}\bigr\Vert _{F}^{2}
\bigr),
\end{eqnarray*}
for some $\delta_{K}^{\prime}$ only depending on $K$ and some absolute
constant $C_{5}$. Since $\Vert U_{i}\Lambda_{t}U_{i}^{T}-V_{0,\bar
{B}}\Lambda
_{0}V_{0,\bar{B}}^{T}\Vert _{F}\geq j\bar{M}\varepsilon$, we have
\begin{eqnarray*}
P_{\Sigma}^{n}\phi_{B,t,j,i}&\leq& 2\exp
\bigl(-C_{5}\delta _{K}^{\prime
}nj^{2}
\bar{M}^{2}\varepsilon^{2} \bigr), \\
\sup_{\Gamma\in
H_{B,t,j,i}}P_{\Gamma}^{n}(1-
\phi_{B,t,j,i})&\leq& 2\exp \bigl(-C_{5}\delta _{K}^{\prime}nj^{2}
\bar{M}^{2}\varepsilon^{2} \bigr).
\end{eqnarray*}

Now we are ready to integrate these little tests step by step for each
index. For each $(B,t,j)$, define
\[
\phi_{B,t,j}=\max_{1\leq i\leq N_{j}}\phi_{B,t,j,i},
\]
and we have
\begin{eqnarray*}
P_{\Sigma}^{n}\phi_{B,t,j} &\leq&\sum
_{i=1}^{N_{j}}P_{\Sigma
}^{n}\phi
_{B,t,j,i}
\\
&\leq&2N_{j}\exp \bigl(-C_{5}\delta_{K}^{\prime}nj^{2}
\bar {M}^{2}\varepsilon ^{2} \bigr)
\\
&\leq&2\exp \biggl(-C_{5}\delta_{K}^{\prime}j^{2}
\bar {M}^{2}n\varepsilon ^{2}+(A+1)r^{2}s\log
\bigl(48\sqrt{2}\delta_{K}^{-1}K \bigr)\\
&&\hspace*{123pt}{}+r^{2}\log
(6%
\sqrt{r} )+\frac{1}{2}r^{2}\log n \biggr).
\end{eqnarray*}
Since we assume $r\vee\log n\leq m\log p$ and $r\leq ms$, we have $%
r^{2}s\leq mn\varepsilon^{2}$,\break $r^{2}\log (6\sqrt{r} )\leq
mn\varepsilon
^{2}$ and $r^{2}\log n\leq m^{2}n\varepsilon^{2}$. Hence
\begin{eqnarray*}
P_{\Sigma}^{n}\phi_{B,t,j} &\leq&2\exp \bigl(-
\bigl(C_{5}\delta _{K}^{\prime}j^{2}
\bar{M}^{2}-(A+1)m\log \bigl(48\sqrt{2}\delta _{K}^{-1}K%
 \bigr)-m-m^{2}/2 \bigr)n\varepsilon^{2} \bigr)
\\
&\leq&2\exp \bigl(-\tfrac{1}{2}C_{5}\delta_{K}^{\prime}j^{2}
\bar{M} 
^{2}n\varepsilon^{2} \bigr),
\end{eqnarray*}
as long as we pick
\[
\bar{M}^{2}\geq2C_{5}^{-1}\delta_{K}^{\prime}{}^{-1}(A+1)m
\log \bigl(48%
\sqrt{2}\delta_{K}^{-1}K
\bigr)+2C_{5}^{-1}\delta_{K}^{\prime
}{}^{-1}m+C_{5}^{-1}
\delta_{K}^{\prime}{}^{-1}m^{2}.
\]
In addition, for each $(B,t,j)$,
\[
\sup_{\Gamma\in H_{1,B,t,j}}P_{\Gamma}^{n}(1-
\phi_{B,t,j})\leq 2\exp \bigl(%
-C_{5}
\delta_{K}^{\prime}j^{2}\bar{M}^{2}n
\varepsilon^{2} \bigr).
\]

For each $(B,t)$, we define
\[
\phi_{B,t}=\max_{j}\phi_{B,t,j},
\]
whose errors are bounded as follows:
\begin{eqnarray*}
P_{\Sigma}^{n}\phi_{B,t} &\leq&\sum
_{j}P_{\Sigma}^{n}\phi_{B,t,j}
\\
&\leq&2\sum_{j}\exp \biggl(-\frac{1}{2}C_{5}
\delta_{K}^{\prime
}j^{2}\bar{M}%
^{2}n
\varepsilon^{2} \biggr)
\\
&\leq&3\exp \biggl(-\frac{1}{2}C_{5}\delta_{K}^{\prime}
\bar{M}%
^{2}n\varepsilon^{2} \biggr)
\end{eqnarray*}
and
\[
\sup_{\Gamma\in H_{B,t}}P_{\Gamma}^{n}(1-
\phi_{B,t})\leq2\exp \bigl(%
-C_{5}
\delta_{K}^{\prime}\bar{M}^{2}n\varepsilon^{2}
\bigr).
\]

For each $B$, we define
\[
\phi_{B}=\max_{1\leq t\leq T}\phi_{B,t},
\]
and we have the errors bounded by
\begin{eqnarray*}
P_{\Sigma}^{n}\phi_{B} &\leq&\sum
_{t=1}^{T}P_{\Sigma}^{n}\phi
_{B,t}
\\
&\leq&3\exp \biggl(-\frac{1}{2}C_{5}\delta_{K}^{\prime}
\bar{M}%
^{2}n\varepsilon^{2}+\log T \biggr)
\\
&\leq&3\exp \biggl(-\frac{1}{2}C_{5}\delta_{K}^{\prime}
\bar{M}%
^{2}n\varepsilon^{2}+r\log\bigl(12 K
\sqrt{r}\varepsilon^{-1}\bigr) \biggr)
\\
&\leq&3\exp \biggl(-\frac{1}{4}C_{5}\delta_{K}^{\prime}
\bar{M}%
^{2}n\varepsilon^{2} \biggr)
\end{eqnarray*}
and
\[
\sup_{\Gamma\in H_{B}}P_{\Gamma}^{n}(1-
\phi_{B})\leq2\exp \bigl(%
-C_{5}
\delta_{K}^{\prime}\bar{M}^{2}n\varepsilon^{2}
\bigr).
\]
Finally, the ultimate test is defined as
\[
\phi=\max_{B}\phi_{B},
\]
with type I error $P_{\Sigma}^{n}\phi$ bounded by
\begin{eqnarray*}
\sum_{B}P_{\Sigma}^{n}
\phi_{B} &\leq& \Biggl( \sum_{q=1}^{[\operatorname{Ars}]}{
\pmatrix{p
\cr
q}} \Biggr) 3\exp \biggl(-\frac{1}{4%
}C_{5}
\delta_{K}^{\prime}\bar{M}^{2}n\varepsilon^{2}
\biggr)
\\
&\leq&3\operatorname{Ars}\exp (\operatorname{Ars}\log p )\exp \biggl(-\frac
{1}{4}C_{5}
\delta _{K}^{\prime}\bar{M}^{2}n\varepsilon^{2}
\biggr)
\\
&\leq&3\exp (2\operatorname{Ars}\log p )\exp \biggl(-\frac{1}{4}C_{5}\delta
_{K}^{\prime}\bar{M}^{2}n\varepsilon^{2}
\biggr)
\\
&\leq&3\exp \biggl(- \biggl(\frac{1}{4}C_{5}
\delta_{K}^{\prime}\bar {M}^{2}-2A%
 \biggr)n\varepsilon^{2} \biggr)
\\
&\leq&3\exp \biggl(-\frac{1}{8}C_{5}\delta_{K}^{\prime}
\bar{M}%
^{2}n\varepsilon^{2} \biggr),
\end{eqnarray*}
as long as we choose $\bar{M}^{2}\geq16\delta_{K}^{\prime
}{}^{-1}C_{5}^{-1}A$, and for type II error we have
\[
\sup_{\Gamma\in\mathcal{H}_{1}}P_{\Gamma}^{n}(1-\phi)\leq2\exp
\bigl(-C_{5}\delta _{K}^{\prime}
\bar{M}^{2}n\varepsilon^{2} \bigr).
\]
Thus the proof is complete.
\end{pf*} 

\section*{Acknowledgments}
We thank the referees and the Associate Editor for giving valuable and
insightful suggestions, which lead to significant improvement of the paper.

\begin{supplement}[id=suppA]
\stitle{Supplement to ``Rate-optimal posterior contraction for sparse PCA''}
\slink[doi]{10.1214/14-AOS1268SUPP} 
\sdatatype{.pdf}
\sfilename{aos1268\_supp.pdf}
\sdescription{In the supplementary text [\citeauthor {gaosupp}
(\citeyear {gaosupp})],
we present proofs of Proposition~\ref{prop:pointestimate}, Lemmas \ref
{lem:denominator},  \ref{lem:specconcen}, \ref{lem:dcover},
Theorem~\ref{thmm:rankone}, Proposition~\ref{prop:chainrule} and Lemma~\ref{lem:testfrob}.}
\end{supplement}


\printaddresses
\end{document}